\documentclass[reqno]{amsart}

\usepackage{hyperref}
\usepackage{amssymb,dsfont}

\newcommand{\intn}[2]{\ensuremath{\{  #1  ,\ldots,  #2 \}}}
\newcommand{\gene}[1]{\ensuremath{\langle #1 \rangle}}

\def\grint{\mathrm{G}}

\def\Id{\mathrm{Id}}

\def\R{\mathbb{R}}
\def\N{\mathbb{N}}
\def\Z{\mathbb{Z}}
\def\T{\mathbb{T}}

\DeclareMathAlphabet{\mathpzc}{OT1}{pzc}{m}{it}
\def\dist{\mathpzc{d}}

\def\hau{\mathcal{H}}

\def\netm{\mathcal{M}}

\def\trib{\mathcal{F}}
\def\tribbis{\mathcal{G}}
\def\prob{\mathbb{P}}
\def\esp{\mathbb{E}}
\def\dd{\mathrm{d}}

\newcommand{\diam}[1]{\ensuremath{| #1 |}}

\def\eps{\varepsilon}

\def\jauge{\mathfrak{D}}

\def\etiq{\mathcal{U}}
\def\raci{\varnothing}
\def\ind{\mathds{1}}

\renewcommand{\theenumi}{\alph{enumi}}

\theoremstyle{plain}
\newtheorem{thm}{Theorem}
\newtheorem{prp}{Proposition}
\newtheorem{lem}[prp]{Lemma}

\theoremstyle{definition}
\newtheorem*{df}{Definition}

\theoremstyle{remark}
\newtheorem{rem}{Remark}
	
\begin{document}

\title[Random wavelet series based on a tree-indexed Markov chain]{Random wavelet series based on\\ a tree-indexed Markov chain}
\author{Arnaud Durand}
\address{Laboratoire d'Analyse et de Math\'ematiques Appliqu\'ees, Universit\'e Paris XII, 61 av. du G\'en\'eral de Gaulle, 94010 Cr\'eteil Cedex, France.}
\email{a.durand@univ-paris12.fr}

\begin{abstract}
We study the global and local regularity properties of random wavelet series whose coefficients exhibit correlations given by a tree-indexed Markov chain. We determine the law of the spectrum of singularities of these series, thereby performing their multifractal analysis. We also show that almost every sample path displays an oscillating singularity at almost every point and that the points at which a sample path has at most a given H\"older exponent form a set with large intersection.
\end{abstract}

\maketitle

\section{Introduction}\label{hmmintro}

Wavelets emerged in the 1980s as a powerful tool for signal processing, see \cite{Grossmann:1984kx,Holschneider:1989uq,Kronland-Martinet:1988fj}. Since then, they have found many applications in this field, such as estimation, detection, classification, compression, filtering and synthesis, see~e.g. \cite{Basseville:1992fj,Chipman:1997fk,Donoho:1995lr,Flandrin:1992kx,Lee:1996qy,Shapiro:1993uq}. In these papers, the coefficients are implicitly assumed to be independent of one another and the exposed methods are based on scalar transformations on each wavelet coefficient of the considered signal. Nonetheless, it was observed that the wavelet coefficients of many real-world signals exhibit some correlations. In particular, the large wavelet coefficients tend to propagate across scales at the same locations, see \cite{Mallat:1992yq,Mallat:1992vn}. In image processing, this phenomenon is related to the fact that the contours of a picture generate discontinuities. Therefore, methods that exploit dependencies between wavelet coefficients should yield better results in the applications.

In order to develop such methods, M. Crouse, R. Nowak and R. Baraniuk \cite{Crouse:1998sh} introduced a simple probabilistic model allowing to capture correlations between wavelet coefficients : the {\em Hidden Markov Tree} (HMT) model, which we now briefly recall. Let us consider a Markov chain $X$ indexed by the binary tree with state space $\{0,1\}$. Basically, the Markov property enjoyed by $X$ means that if the state $X_{u}$ of a vertex $u$ is $x\in\{0,1\}$, then the states of the two sons of $u$ are chosen independently according to the transition probabilities from $x$. Conditionally on the Markov chain $X$, the wavelet coefficients are independent and the wavelet coefficient indexed by a given dyadic interval $\lambda$ is a centered Gaussian random variable whose variance is large (resp.~small) if the state of the vertex of the binary tree that corresponds to $\lambda$ is $1$ (resp.~$0$). The correlations between the wavelet coefficients are thus given by the underlying Markov chain. Furthermore, the unconditional law of each coefficient is a Gaussian mixture. This property agrees with the observation that the wavelet coefficient histogram of a real-world signal is usually more peaky at zero and heavy-tailed than the Gaussian. Moreover, let us mention that the HMT model was used in image processing, see \cite{Choi:2001ys,Diligenti:2001rt} for instance.

In this paper, we investigate the pointwise regularity properties of the sample paths of a model of random wavelet series which is closely related to the HMT model. The regularity of a function at a given point is measured by its H\"older exponent, which is defined as follows, see \cite{Jaffard:2004fh}.

\begin{df}[H\"older exponent]
Let $f$ be a function defined on $\R$ and let $x\in\R$. The H\"older exponent $h_{f}(x)$ of $f$ at $x$ is the supremum of all $h>0$ such that there are two reals $c>0$ and $\delta>0$ and a polynomial $P$ enjoying
\[
\forall x'\in [x-\delta,x+\delta] \qquad |f(x')-P(x'-x)|\leq c\,|x'-x|^h.
\]
\end{df}

Since we are interested in local properties, it is more convenient to work with wavelets on the torus $\T=\R/\Z$, so that the random wavelet series that we study throughout the paper is in fact a random process on $\T$. This process is defined in Section \ref{presentation} and is denoted by $R$. Let $\phi$ be the canonical surjection from $\R$ to $\T$. Note that $R\circ\phi$ is a one-periodic random function defined on $\R$. The H\"older exponent of $R$ at a point $x\in\T$ is then defined by $h_{R}(x)=h_{R\circ\phi}(\dot x)$, where $\dot x$ is a real number such that $\phi(\dot x)=x$. Equivalently, $h_{R}(x)$ is the supremum of all $h>0$ such that
\[
\dist(x',x)\leq\delta \qquad\Longrightarrow\qquad |R(x')-P(x'-x)|\leq c\,\dist(x',x)^h
\]
for all $x'\in\T$, some positive reals $c$ and $\delta$ and some function $P$ on $\T$ such that $P\circ\phi$ coincides with a polynomial in a neighborhood of zero, where $\dist$ denotes the quotient distance on the torus $\T$.

Our main purpose is to perform the {\em multifractal analysis} of the sample paths of the random wavelet series $R$. This amounts to studying the size properties of the {\em iso-H\"older set}
\begin{equation}\label{defEhhmm}
E_{h}=\left\{ x\in\T \:|\: h_{R}(x)=h \right\}
\end{equation}
for every $h\in [0,\infty]$. More precisely, in Subsection \ref{lawspecsubsec}, we give the law of the mapping $d_{R}:h\mapsto\dim E_{h}$, where $\dim$ stands for Hausdorff dimension. This mapping is called the {\em spectrum of singularities} of the process $R$.

A remarkable property, due to the correlations between wavelet coefficients, is that the spectrum of singularities of the process $R$ is itself a random function. None of the multifractal stochastic processes studied up to now, such as the usual L\'evy processes \cite{Jaffard:1999fg}, the L\'evy processes in multifractal time \cite{Barral:2005ef} or the random wavelet series with independent coefficients introduced by J.-M. Aubry and S. Jaffard \cite{Aubry:2002up,Jaffard:2000mk}, enjoys this property. The random wavelet series based on multifractal measures studied by J. Barral and S. Seuret \cite{Barral:2005fr} do not satisfy it either, even though their wavelet coefficients exhibit strong correlations. This remark also holds for the wavelet series based on branching processes introduced by A. Brouste \cite{Brouste:2006lr}. In this last model, for any dyadic interval $\lambda$, the wavelet coefficient indexed by $\lambda$ is either Gaussian or zero depending on whether or not the vertex of the binary tree that corresponds to $\lambda$ belongs to a certain Galton-Watson tree. This model can be seen as a particular case of the HMT model by assuming that the underlying Markov chain cannot map to the state $1$ a vertex whose father is mapped to the state $0$. Let us also mention that A. Brouste, in collaboration with geophysicists, used this model to study the surface roughness of certain rocks, see \cite{Brouste:2007rt}.

For certain values of the parameters of the model, the sample paths of the random wavelet series $R$ enjoy the remarkable property that almost every of them displays an {\em oscillating singularity} at almost every point of the torus.  We refer to Subsection \ref{oscsingsubsec} for details. Note that this property also holds for the models of random wavelet series studied by J.-M. Aubry and S. Jaffard \cite{Aubry:2002up,Jaffard:2000mk}. Conversely, the random wavelet series considered by J. Barral and S. Seuret \cite{Barral:2005fr} and those introduced by A. Brouste \cite{Brouste:2006lr} do not exhibit any oscillating singularity.

Certain random sets related to the iso-H\"older sets $E_{h}$ enjoy a notable geometric property which was introduced by K. Falconer \cite{Falconer:1994hx}. To be specific, we establish that, for certain values of $h$, the sets
\begin{equation}\label{deftildeEhhmm}
\widetilde E_{h}=\left\{ x\in\T \:\bigl|\: h_{R}(x)\leq h \right\}
\end{equation}
are almost surely {\em sets with large intersection}, see Subsection \ref{largeint}. In particular, this implies that they are locally everywhere of the same size, in the sense that the Hausdorff dimension of the set $\widetilde E_{h}\cap V$ does not depend on the choice of the nonempty open subset $V$ of the torus. This also implies that the size properties of the sets $\widetilde E_{h}$ are not altered by taking countable intersections. In fact, the Hausdorff dimension of the intersection of countably many sets with large intersection is equal to the infimum of their Hausdorff dimensions. This property is somewhat counterintuitive in view of the fact that the intersection of two subsets of the torus with Hausdorff dimensions $s_{1}$ and $s_{2}$ respectively is usually expected to be $s_{1}+s_{2}-1$ (see \cite[Chapter~8]{Falconer:2003oj} for precise statements). The occurrence of sets with large intersection in the theory of Diophantine approximation and that of dynamical systems was pointed out by several authors, see \cite{Durand:2007uq,Durand:2006uq,Falconer:1994hx} and the references therein. Their use in multifractal analysis of stochastic processes is more novel and was introduced by J.-M. Aubry and S. Jaffard. Indeed, they established in \cite{Aubry:2002up} that sets with large intersection arise in the study of the H\"older singularity sets of certain random wavelet series. As shown in \cite{Durand:2007fk}, such sets also appear in the study of the singularity sets of L\'evy processes.

The rest of the paper is organized as follows. In Section \ref{presentation} we give a precise definition of the model of random wavelet series that we study. Our main results are stated in Section \ref{results} and are proven in Sections \ref{hmmprel} to \ref{oscsing}.

\section{Presentation of the process}\label{presentation}

In order to define the process that we study, let us introduce some notations. Throughout the paper, $\N$ (resp.~$\N_{0}$) denotes the set of positive (resp.~nonnegative) integers and $\Lambda$ is the collection of the dyadic intervals of the torus $\T$, that is, the sets of the form $\lambda=\phi(2^{-j}(k+[0,1)))$ with $j\in\N_{0}$ and $k\in\intn{0}{2^j-1}$. The integer $\gene{\lambda}=j$ is called the generation of $\lambda$. Furthermore, let us consider a wavelet $\psi$ in the Schwartz class (see \cite{Lemarie-Rieusset:1986lr}). For any dyadic interval $\lambda=\phi(2^{-j}(k+[0,1)))\in\Lambda$, let $\Psi_{\lambda}$ denote the function on $\T$ which corresponds to the one-periodic function
\[
x\mapsto\sum_{m\in\Z} \psi(2^j(x-m)-k).
\]
Then, the functions $2^{\gene{\lambda}/2}\Psi_{\lambda}$, together with the constant function equal to one on $\T$, form an orthonormal basis of $L^2(\T)$, see \cite{Meyer:1990qj}.

Recall that there is a one-to-one correspondence between the set $\Lambda$ of all dyadic intervals of the torus and the set
\[
\etiq = \{\raci\}\cup\bigcup_{j=1}^\infty \{0,1\}^j.
\]
The set $\etiq$ is formed of the empty word $\raci$ and the words $u=u_{1}\ldots u_{j}$ of finite length $j\geq 1$ in the alphabet $\{0,1\}$. The integer $\gene{u}=j$ is called the generation of $u$. In addition, let $\gene{\raci}=0$ and $\etiq^*=\etiq\setminus\{\raci\}$. For every word $u\in\etiq^*$, the word $\pi(u)=u_{1}\ldots u_{\gene{u}-1}$ is called the father of $u$. Then, the directed graph with vertex set $\etiq$ and with arcs $(\pi(u),u)$, for $u\in\etiq^*$, is a binary tree rooted at $\raci$. To be specific, the bijection from $\etiq$ to $\Lambda$ is
\[
u \mapsto \lambda_{u}=\phi\left(\sum_{j=1}^{\gene{u}} u_{j} 2^{-j}+[0,2^{-\gene{u}})\right).
\]
Thus, for every $\lambda\in\Lambda$, there is a unique vertex $u_{\lambda}\in\etiq$ such that $\lambda_{u_{\lambda}}=\lambda$.

In the following, $X=(X_{u})_{u\in\etiq}$ denotes a $\{0,1\}$-valued stochastic process indexed by the binary tree $\etiq$. The $\sigma$-field
\[
\tribbis_{u}=\sigma(X_{v},\,v\in\etiq\setminus(u\etiq^*))
\]
can be considered as the {\em past} before $u$, because the set $u\etiq^*$ is composed of all the {\em descendants} of the vertex $u$ in the binary tree $\etiq$, that is, words of the form $uv$ with $v\in\etiq^*$. Conversely, the {\em future} after $u$ begins with its two sons $u0$ and $u1$. For any integer $j\geq 0$, let $\nu_{0,j}$ and $\nu_{1,j}$ denote two probability measures on $\{0,1\}^2$. From now on, we assume that the process $X$ is a Markov chain with transition probabilities given by the measures $\nu_{0,j}$ and $\nu_{1,j}$. This means that the following {\em Markov condition} holds:
\renewcommand{\theenumi}{\Alph{enumi}}
\begin{enumerate}\setcounter{enumi}{12}
\item\label{condMarkov} For any vertex $u\in\etiq$ and any subset $A$ of $\{0,1\}^2$,
\[
\prob((X_{u0},X_{u1})\in A \:|\: \tribbis_{u})=\nu_{X_{u},\gene{u}}(A).
\]
\end{enumerate}
\renewcommand{\theenumi}{\alph{enumi}}
Equivalently, the conditional distribution of the vector $(X_{u0},X_{u1})$ conditionally on the past $\tribbis_{u}$ is the probability measure $\nu_{X_{u},\gene{u}}$. Thus, the Markov condition satisfied by $X$ is conceptually similar to that enjoyed by inhomogeneous discrete time Markov chains.

The random wavelet series $R$ that we study is then defined by
\begin{equation*}\begin{split}
R &=\sum_{u\in\etiq} \left( 2^{-\underline{h}\,\gene{u}}\ind_{\{X_{u}=1\}} + 2^{-\overline{h}\,\gene{u}}\ind_{\{X_{u}=0\}} \right)\Psi_{\lambda_{u}} \\
&=\sum_{\lambda\in\Lambda} \left( 2^{-\underline{h}\,\gene{\lambda}}\ind_{\{X_{u_{\lambda}}=1\}} + 2^{-\overline{h}\,\gene{\lambda}}\ind_{\{X_{u_{\lambda}}=0\}} \right)\Psi_{\lambda},
\end{split}\end{equation*}
where $0<\underline{h}<\overline{h}\leq\infty$. The mean value of every sample path of $R$ vanishes. Moreover, for any $\lambda\in\Lambda$, the wavelet coefficient indexed by $\lambda$ is
\begin{equation}\label{defClambda}
\begin{cases}
C_{\lambda}=2^{-\underline{h}\,\gene{\lambda}} & \text{if }X_{u_{\lambda}}=1 \\[2mm]
C_{\lambda}=2^{-\overline{h}\,\gene{\lambda}} & \text{if }X_{u_{\lambda}}=0.
\end{cases}
\end{equation}
The coefficient $C_{\lambda}$ should be considered as large in the first case and small in the second one (it even vanishes if $\overline{h}=\infty$ and $\gene{\lambda}\geq 1$). The model that we study is thus defined in the same way as the HMT model, except that, conditionally on the underlying Markov chain, the wavelet coefficients are deterministic instead of Gaussian.

Note that the Markov condition (\ref{condMarkov}) implies that the probability measures $\nu_{1,j}$ affect the propagation of the large wavelet coefficients of $R$ across scales, while the probability measures $\nu_{0,j}$ govern their appearance. In what follows, the influence of each of these two phenomena is reflected by the values of the parameters
\begin{equation}\label{defgammaj}
\gamma_{j}=2\,\nu_{1,j}(\{(1,1)\})+\nu_{1,j}(\{(1,0),(0,1)\})
\end{equation}
and
\begin{equation}\label{defetaj}
\eta_{j}=1-\nu_{0,j}(\{(0,0)\}),
\end{equation}
respectively. Indeed, for any integer $j\geq 0$, $\gamma_{j}$ is the expected number of sons mapped to the state $1$ of a vertex with generation $j$ that is mapped to the state $1$ and $\eta_{j}$ is the probability that a vertex with generation $j$ that is mapped to the state $0$ has at least one son mapped to the state $1$.

\section{Statement of the results}\label{results}

\subsection{Preliminary remarks}

As shown by (\ref{defClambda}), for every dyadic interval $\lambda\in\Lambda$, the modulus of the wavelet coefficient $C_{\lambda}$ of $R$ is at most $2^{-\underline{h}\gene{\lambda}}$. A standard result of \cite{Meyer:1990qj} then implies that $R$ belongs to the H\"older space of uniform regularity $C^{\underline{h}}(\T)$. It follows that the H\"older exponent $h_{R}(x)$ defined in Section \ref{hmmintro} is at least $\underline{h}$ for every point $x\in\T$. As shown by Theorems \ref{bifrac} and \ref{multifrac} below, the H\"older exponent of $R$ is actually greater than $\underline{h}$ at most locations.

The H\"older exponent of $R$ at any point of the torus is highly related to the size of the wavelet coefficients which are indexed by the dyadic intervals located around $x$. More precisely, it can be computed using the following proposition, which is a straightforward consequence of \cite[Proposition 1.3]{Jaffard:1996kt}.

\begin{prp}\label{critponclem}
\noindent Let $h\in (0,\infty)$ and $x\in\T$.
\begin{enumerate}
\item\label{critponclem1} If $h_{R}(x)>h$, then
\begin{equation}\label{critponc}
\exists \kappa>0 \quad \forall \lambda\in\Lambda \qquad |C_{\lambda}| \leq \kappa\,(2^{-\gene{\lambda}}+\dist(x,x_{\lambda}))^h.
\end{equation}
\item If (\ref{critponc}) holds, then $h_{R}(x)\geq h$.
\end{enumerate}
\end{prp}

By virtue of (\ref{defClambda}), for each dyadic interval $\lambda\in\Lambda$, the modulus of the wavelet coefficient $C_{\lambda}$ of $R$ is at least $2^{-\overline{h}\gene{\lambda}}$. According to Proposition \ref{critponclem}(\ref{critponclem1}), the H\"older exponent of $R$ is less than $\overline{h}$ everywhere. It follows that
\begin{equation}\label{underhleqhRxleqoverh}
\forall x\in\T \qquad \underline{h}\leq h_{R}(x)\leq\overline{h}.
\end{equation}
As a consequence, for every sample path of the random wavelet series $R$ and every $h\in [0,\underline{h})\cup(\overline{h},\infty]$, the iso-H\"older set $E_{h}$ defined by (\ref{defEhhmm}) is empty. Conversely, for each $h\in [\underline{h},\overline{h}]$, the set $E_{h}$ need not be empty and its size properties are described by Theorems \ref{bifrac} and \ref{multifrac} below.

Before stating these theorems, we need to introduce some further notations related to the parameters $\gamma_{j}$ defined by (\ref{defgammaj}) and governing the propagation of the large wavelet coefficients of the process $R$. To be specific, let
\[
\underline{j}=\inf\{ j_{0}\geq 0 \:|\: \forall j\geq j_{0} \quad \gamma_{j}>0 \}
\]
and
\[
\begin{cases}
\theta=\liminf\limits_{j\to\infty} \frac{\log\gamma_{\underline{j}}+\ldots+\log\gamma_{j}}{j\log 2} & \text{if }\underline{j}<\infty \\[2mm]
\theta=-\infty & \text{if }\underline{j}=\infty.
\end{cases}
\]
Note that $\theta$ is at most one. Moreover, $\gamma_{j}$ being the expected number of sons mapped to the state $1$ of a vertex with generation $j$ that is mapped to the state $1$, the parameter $\theta$ expresses the trend with which large wavelet coefficients propagate across scales. For any integer $j\geq 0$, we also need to consider the number
\[
\varsigma_{j}=2\sum_{n=j}^\infty \frac{\nu_{1,n}(\{1,1\})}{\gamma_{n}\prod\limits_{\ell=j}^n \gamma_{\ell}},
\]
which naturally occurs in the study of a specific random fractal set related to the process $R$, see Lemma \ref{prprecaptheta} below.

Throughout the rest of the paper, we suppose that $\theta$ is less than one. This assumption implies that the large wavelet coefficients of $R$ cannot propagate too much across scales. It is not very restrictive, in view of the fact that the decomposition of a typical real-world signal in a wavelet basis has very few large coefficients.

\subsection{Law of the spectrum of singularities}\label{lawspecsubsec}

Theorems \ref{bifrac} and \ref{multifrac} below give the law of the spectrum of singularities of the sample paths of the random wavelet series $R$. Recall that it is the mapping $d_{R}:h\mapsto\dim E_{h}$, where $E_{h}$ is defined by (\ref{defEhhmm}) and $\dim$ stands for Hausdorff dimension.

With a view to recalling the definition of Hausdorff dimension, let us first define the notion of Hausdorff measure on the torus. To this end, let $\jauge$ denote the set of all nondecreasing functions $g$ defined on a neighborhood of zero and enjoying $\lim_{0^+}g=g(0)=0$. Any function in $\jauge$ is called a gauge function. For every $g\in\jauge$, the Hausdorff $g$-measure of a subset $F$ of $\T$ is defined by
\[
\hau^g(F)=\lim_{\eps\downarrow 0}\uparrow \hau^g_{\eps}(F) \qquad\text{with}\qquad \hau^g_{\eps}(F)=\inf_{F\subseteq\bigcup_{p} U_{p}\atop\diam{U_{p}}<\eps} \sum_{p=1}^\infty g(\diam{U_{p}}).
\]
The infimum is taken over all sequences $(U_{p})_{p\geq 1}$ of sets with $F\subseteq\bigcup_{p} U_{p}$ and $\diam{U_{p}}<\eps$ for all $p\geq 1$, where $\diam{\cdot}$ denotes diameter (in the sense of the quotient distance on the torus). Note that $\hau^g$ is a Borel measure on $\T$, see \cite{Rogers:1970wb}. The Hausdorff dimension of a nonempty set $F\subseteq\T$ is then defined by
\[
\dim F=\sup\{s\in (0,1) \:|\: \hau^{\Id^s}(F)=\infty\}=\inf\{s\in (0,1) \:|\: \hau^{\Id^s}(F)=0\},
\]
where $\Id^s$ denotes the function $r\mapsto r^s$ and with the convention that $\sup\emptyset=0$ and $\inf\emptyset=1$, see \cite{Falconer:2003oj}. In addition, we agree with the usual convention that the empty set $\emptyset$ has Hausdorff dimension $-\infty$.

The pointwise regularity properties of the sample paths of the process $R$ depend significantly on whether or not the series $\sum_{j} 2^j \eta_{j}$ converges, where $\eta_{j}$ is given by (\ref{defetaj}). If it converges, then the probabilities with which large wavelet coefficients can appear at a given scale are quite low, so that the sample paths of $R$ are regular in their irregularity, as shown by the following result. In its statement, the function $\Phi_{j}$ is the generating function of the cardinality of the set
\begin{equation}\label{defSj}
S_{j}=\left\{ u\in\etiq \:|\: \gene{u}=j \text{ and } X_{u}=1 \right\}.
\end{equation}
Equivalently, $\Phi_{j}(z)=\esp[z^{\# S_{j}}]$ for any complex number $z$ and any integer $j\geq 0$. The functions $\Phi_{j}$ may be calculated recursively in terms of the transition probabilities of the Markov chain $X$. Indeed, the Markov condition (\ref{condMarkov}) implies that for any integer $j\geq 0$,
\[
\esp[z^{\# S_{j+1}}\:|\:\trib_{j}]=\left(\int_{\{0,1\}^2} z^{x_{0}+x_{1}}\,\nu_{1,j}(\dd x)\right)^{\# S_{j}}\left(\int_{\{0,1\}^2} z^{x_{0}+x_{1}}\,\nu_{0,j}(\dd x)\right)^{2^j-\# S_{j}},
\]
where $\trib_{j}$ is the $\sigma$-field generated by the variables $X_{u}$, for $u\in\etiq$ with $\gene{u}\leq j$.

\begin{thm}\label{bifrac}
Let us suppose that $\sum_{j} 2^j \eta_{j}<\infty$ and $\theta<1$.
\begin{enumerate}
\item If $\theta<0$, then with probability one, $h_{R}(x)=\overline{h}$ for all $x\in\T$. Therefore, with probability one, for all $h\in [0,\infty]$,
\[
\begin{cases}
d_{R}(h)=1 & \text{if }h=\overline{h} \\[2mm]
d_{R}(h)=-\infty & \text{else}.
\end{cases}
\]
\item If $\theta\geq 0$, then with probability one, for all $h\in [0,\infty]$,
\[
\begin{cases}
d_{R}(h)\in\{-\infty,\theta\} & \text{if }h=\underline{h} \\[2mm]
d_{R}(h)=1 & \text{if }h=\overline{h} \\[2mm]
d_{R}(h)=-\infty & \text{else}.
\end{cases}
\]
Moreover, if there is an integer $j_{*}\geq 0$ such that $\nu_{1,j}(\{(0,0)\})=0$ for any $j\geq j_{*}$, then
\[
\prob(d_{R}(\underline{h})=-\infty)=\Phi_{j_{*}}(0)\cdot\prod_{j=j_{*}}^\infty (1-\eta_{j})^{2^j}.
\]
If not, $\prob(d_{R}(\underline{h})=-\infty)$ is positive and it is equal to one if and only if
\[
\varsigma_{\underline{j}}=\infty
\qquad\text{or}\qquad
\liminf_{j\to\infty}\prod_{\ell=\underline{j}}^j \gamma_{\ell}=0
\qquad\text{or}\qquad
\left\{\begin{array}{l}
\Phi_{\underline{j}}(0)=1 \\[2mm]
\forall j\geq\underline{j} \quad \eta_{j}=0.
\end{array}\right.
\]
\end{enumerate}
\end{thm}

\begin{rem}
As $\theta<1$, Theorem \ref{bifrac} implies that the iso-H\"older set $E_{\overline{h}}$ has full Lebesgue measure in the torus with probability one. Thus, the H\"older exponent of almost every sample path of the process $R$ is $\overline{h}$ almost everywhere.
\end{rem}

\begin{rem}
The process $R$ enables to provide a partial answer to a question raised by S. Jaffard in \cite{Jaffard:2004hy}. For many examples of random processes $F$, such as the usual L\'evy processes \cite{Jaffard:1999fg}, the L\'evy processes in multifractal time \cite{Barral:2005ef} or several models of random wavelet series \cite{Aubry:2002up,Jaffard:2000mk}, even though the function $x \mapsto h_{F}(x)$ is random, the spectrum of singularities of $F$ is a deterministic function. Of course, this property does not hold in general: consider for instance a fractional Brownian motion whose Hurst parameter follows a Bernoulli law. The random wavelet series $R$ supplies a more elaborate example. Indeed, Theorem \ref{bifrac} shows that its spectrum of singularities may be random. Theorem \ref{multifrac} below indicates that this property still holds if the series $\sum_{j} 2^j \eta_{j}$ diverges.
\end{rem}

If the series $\sum_{j} 2^j \eta_{j}$ diverges, large wavelet coefficients can appear with a relatively large probability at each scale, which makes the sample paths of $R$ very irregular in their irregularity, as shown by the following result. In its statement, $\widetilde h$ is defined by
\[
\widetilde h=\inf\left\{ h>0 \:\biggl|\: \sum\nolimits_{j} 2^{(1-\underline{h}/h)j}\eta_{j}=\infty \right\},
\]
which is clearly greater than or equal to $\underline{h}$.

\eject

\begin{thm}\label{multifrac}
Let us suppose that $\sum_{j} 2^j \eta_{j}=\infty$ and $\theta<1$.
\begin{enumerate}
\item If $\widetilde h<\overline{h}$, then with probability one, for all $h\in [0,\infty]$,
\[
\begin{cases}
d_{R}(h)=h/\widetilde h & \text{if }\underline{h}<h\leq\widetilde h \\[2mm]
d_{R}(h)=-\infty & \text{if }h<\underline{h}\text{ or }h>\widetilde h.
\end{cases}
\]
\item If $\widetilde h\geq\overline{h}$, then with probability one, for all $h\in [0,\infty]$,
\[
\begin{cases}
d_{R}(h)=h/\widetilde h & \text{if } \underline{h}<h<\overline{h} \\[2mm]
d_{R}(h)=1 & \text{if }h=\overline{h} \\
d_{R}(h)=-\infty & \text{if }h<\underline{h}\text{ or }h>\overline{h}.
\end{cases}
\]
\item If $\underline{h}/\widetilde h\geq\theta$, then with probability one,
\[
d_{R}(\underline{h})=\underline{h}/\widetilde h.
\]
\item If $\underline{h}/\widetilde h<\theta$, then with probability one,
\[
d_{R}(\underline{h})\in\{\underline{h}/\widetilde h,\theta\}.
\]
Moreover, $\prob(d_{R}(\underline{h})=\theta)$ is positive and it is equal to one if and only if either $\sum_{j}2^j\eta_{j}/\varsigma_{j+1}=\infty$ or $\nu_{1,j}(\{(0,0)\})=0$ for all $j$ large enough.
\end{enumerate}
\end{thm}

\begin{rem}
As $\theta<1$, Theorem \ref{multifrac} implies that the iso-H\"older set $E_{\min(\widetilde h,\overline{h})}$ has full Lebesgue measure in the torus with probability one. So, the H\"older exponent of almost every sample path of the random wavelet series $R$ is $\min(\widetilde h,\overline{h})$ almost everywhere. In addition, if $\widetilde h=\underline{h}$, then the set $E_{\widetilde h}$ is almost surely equal to the whole torus, so that the H\"older exponent of $R$ is almost surely $\widetilde h$ everywhere.
\end{rem}

\begin{rem}
For some values of the parameters, the spectrum of singularities of the random wavelet series $R$ need not be concave. Therefore, this spectrum cannot be determined using the multifractal formalisms derived in Besov or oscillation spaces. We refer to \cite{Jaffard:2004fh} for details concerning these multifractal formalisms. Moreover, in general, the spectrum of singularities of $R$ does not coincide with its large deviation spectrum (see e.g. \cite{Levy-Vehel:1997lr} for a fuller exposition), that is, the mapping
\[
h\mapsto \lim_{\eps\downarrow 0}\downarrow\limsup_{j\to\infty}\frac{1}{j}\log_{2}\#\{\lambda\in\Lambda \:|\: \gene{\lambda}=j \text{ and } 2^{-(h+\eps)j}\leq |C_{\lambda}|\leq 2^{-(h-\eps)j}\}.
\]
Indeed, this last function clearly maps any real $h\notin\{\underline{h},\overline{h}\}$ to $-\infty$.
\end{rem}

\begin{rem}
Recall that, owing to Theorem \ref{bifrac}, the spectrum of singularities of $R$ may be random when $\sum_{j} 2^j \eta_{j}<\infty$. Theorem \ref{multifrac} shows that this property still holds when $\sum_{j} 2^j \eta_{j}=\infty$. Specifically, the spectrum of singularities of $R$ is random if and only if $\underline{h}/\widetilde h$ is less than $\theta$, the sum $\sum_{j}2^j\eta_{j}/\varsigma_{j+1}$ is finite and $\nu_{1,j}(\{(0,0)\})$ is positive for infinitely many integers $j\geq 0$.

Let us give an explicit example of probability measures $\nu_{0,j}$ and $\nu_{1,j}$ for which all these conditions hold. Given a real $a\in (0,1)$ and an integer $b\geq 2$, let $p_{0}=2^{-a}$ and, for any integer $j\geq 1$, let $p_{j} = 2^{-a (b^{\lfloor \log_{b} (j+1) \rfloor} - b^{\lfloor \log_{b} j \rfloor})}$, where $\lfloor\,\cdot\,\rfloor$ denotes the floor function. Next, let us consider that the measures $\nu_{1,j}$ are the products
\[
\nu_{1,j}=\left( p_{j}\delta_{1}+(1-p_{j})\delta_{0} \right)^{\otimes 2},
\]
where $\delta_{0}$ and $\delta_{1}$ are the point masses at zero and one, respectively. Clearly, $\nu_{1,j}(\{(0,0)\})$ is positive for all $j$, the number $\underline{j}$ vanishes and $\theta=1-a<1$. For every integer $n\geq 0$, let $j_{n}=b^{n+1}-1$. Observe that, for all $n$ greater than some $n_{0}$, the sum $\varsigma_{j_{n}}$ is at least $2^{a(1-1/b)(j_{n}+1)-1}$ and there exists a real $q_{j_{n}-1}$ such that
\[
2^{-(j_{n}-1)} \leq q_{j_{n}-1} \leq \frac{1}{{j_{n}}^2} \, 2^{\left(a(1-1/b)-1\right) j_{n}}.
\]
Furthermore, let $q_{j-1}=0$ for every integer $j\geq 1$ that is not of the form $j_{n}$ with $n>n_{0}$. Then, let us consider that the probability measures $\nu_{0,j}$ are given by
\[
\nu_{0,j}=\left( q_{j}\delta_{1}+(1-q_{j})\delta_{0} \right)^{\otimes 2}.
\]
It is easy to check that the sum $\sum_{j} 2^j \eta_{j}$ diverges and that $\widetilde h$ is at least $\underline{h}/(a(1-1/b))$. Let us suppose that $a<1/(2-1/b)$. This assumption ensures that $\underline{h}/\widetilde h$ is less than $\theta$. Moreover,
\[
\sum_{n=n_{0}+1}^\infty \frac{2^{j_{n}-1} \eta_{j_{n}-1}}{\varsigma_{j_{n}}} \leq \sum_{n=n_{0}+1}^\infty \frac{2^{a(1-1/b)j_{n}}}{{j_{n}}^2\, 2^{a(1-1/b)(j_{n}+1)-1}} \leq 2^{1-a(1-1/b)}\sum_{j=1}^\infty \frac{1}{j^2}
\]
which ensures the finiteness of the sum $\sum_{j}2^j\eta_{j}/\varsigma_{j+1}$. Theorem \ref{multifrac} finally implies that the spectrum of singularities of the process $R$ is random when the probability measures $\nu_{0,j}$ and $\nu_{1,j}$ are chosen as above.
\end{rem}

\subsection{Oscillating singularities}\label{oscsingsubsec}

Theorem \ref{bifrac} and \ref{multifrac} above give the law of the Hausdorff dimension of the iso-H\"older sets $E_{h}$ defined by (\ref{defEhhmm}). Each set $E_{h}$ is formed of the points at which the H\"older exponent of the process $R$ is $h$. It is possible to provide a more precise description of the pointwise regularity properties of these points. Indeed, a given H\"older exponent $h$ at a point $x\in\T$ can result from many possible local behaviors near $x$. For example, if $h$ is not an even integer, then the {\em cusp} $x'\mapsto |x'-x|^h$ and the {\em chirp}
\begin{equation}\label{chirp}
x'\mapsto |x'-x|^h\sin\frac{1}{|x'-x|^\beta}
\end{equation}
both have H\"older exponent $h$ at $x$, in spite of the fact that their oscillatory behavior is completely different, see \cite{Jaffard:1996kt}.

The oscillating singularity exponent was introduced in \cite{Arneodo:1998ai} in order to describe the oscillatory behavior of a function near a given point and thus to determine if a function behaves rather like a cusp or like a chirp in a neighborhood of a point. It is defined using primitives of fractional order. To be specific, for any $t>0$, any locally bounded function $f$ defined on $\R$ and any $x\in\R$ with $h_{f}(x)<\infty$, let $h^t_{f}(x)$ denote the H\"older exponent at $x$ of the function $(\operatorname{Id}-\Delta)^{-t/2} (\chi f)$, where $\chi$ is a compactly supported smooth function which is equal to one in a neighborhood of $x$ and $(\operatorname{Id}-\Delta)^{-t/2}$ is the operator that corresponds to multiplying by $\xi\mapsto (1+\xi^2)^{-t/2}$ in the Fourier domain.

\begin{df}[oscillating singularity exponent]
Let $f$ be a locally bounded function defined on $\R$ and let $x\in\R$ with $h_{f}(x)<\infty$. The oscillating singularity exponent of $f$ at $x$ is
\[
\beta_{f}(x) = \left.\frac{\partial h^{t}_{f}(x)}{\partial t} \right|_{t=0^+} - 1\in [0,\infty].
\]
If $\beta_{f}(x)>0$, then $f$ is said to display an oscillating singularity at $x$.
\end{df}

It is proven in \cite{Jaffard:1996kt} that if $f$ is defined by (\ref{chirp}), then $h^t_{f}(x)=h+t(\beta+1)$, so that $\beta_{f}(x)=\beta$. As required, the oscillating singularity exponent enables to recover the parameter $\beta$ which governs the oscillatory behavior of a chirp. Note that this exponent is not defined for points at which the H\"older exponent is infinite.

The oscillating singularity exponent $\beta_{R}(x)$ of the random wavelet series $R$ at any point $x\in\T$ such that $h_{R}(x)<\infty$ is then defined in the natural way, that is, $\beta_{R}(x)$ is set to be equal to $\beta_{R\circ\phi}(\dot x)$ for any real number $\dot x$ enjoying $\phi(\dot x)=x$, where $\phi$ is the canonical surjection from $\R$ to $\T$. The following result, which is proven in Section \ref{oscsing}, gives the value of the oscillating singularity exponent of $R$ at every point of the iso-H\"older set $E_{h}$.

\begin{prp}\label{prpexposchmm}
For every $h\in [\underline{h},\overline{h}]$ and every $x\in E_{h}$,
\[
\beta_{R}(x)=\begin{cases}
h/\underline{h}-1 & \text{if }h<\overline{h} \\[2mm]
0 & \text{if }h=\overline{h}<\infty.
\end{cases}
\]
\end{prp}

\begin{rem}
Proposition \ref{prpexposchmm} ensures that the random wavelet series $R$ displays an oscillating singularity at every point of the set $E_{h}$, for any $h\in (\underline{h},\overline{h})$. Moreover, it is necessary to assume the finiteness of $\overline{h}$ for $h=\overline{h}$ in the statement of Proposition \ref{prpexposchmm} because the oscillating singularity exponent is not defined for points at which the H\"older exponent is infinite.
\end{rem}

\begin{rem}
In the case where $\widetilde h<\overline{h}$, Theorem \ref{multifrac} and Proposition \ref{prpexposchmm} ensure that almost every sample path of the random wavelet series $R$ displays an oscillating singularity at almost every point of the torus. Note that this remarkable property is also verified by the models of random wavelet series with independent coefficients which were studied in \cite{Aubry:2002up,Jaffard:2000mk}.
\end{rem}

\subsection{Large intersection properties of the singularity sets}\label{largeint}

For certain values of $h$, the sets $\widetilde E_{h}$ defined by (\ref{deftildeEhhmm}) are sets with large intersection, in the sense that they belong to specific classes $\grint^g(\T)$ of subsets of the torus. These classes are the transposition into the toric setting of the classes $\grint^g(\R)$ of subsets of $\R$ which were introduced in \cite{Durand:2007uq} in order to generalize the original classes of sets with large intersection of K. Falconer \cite{Falconer:1994hx}.

Let us first recall the definition and the basic properties of the classes $\grint^g(\R)$. To begin with, they are defined for functions $g$ in a set denoted by $\jauge_{1}$. This is the set of all gauge functions $g\in\jauge$ such that $r\mapsto g(r)/r$ is positive and nonincreasing on a neighborhood of zero. For any $g\in\jauge_{1}$, let $\eps_{g}$ denote the supremum of all $\eps\in (0,1]$ such that $g$ is nondecreasing on $[0,\eps]$ and $r\mapsto g(r)/r$ is nonincreasing on $(0,\eps]$ and let $\Lambda_{g}$ denote the set of all dyadic intervals of diameter less than $\eps_{g}$, that is, sets of the form $\lambda=2^{-j}(k+[0,1))$ with $j\in\N_{0}$, $k\in\Z$ and $\diam{\lambda}<\eps_{g}$. The outer net measure associated with $g$ is defined by
\[
\forall F\subseteq\R \qquad \netm^g_{\infty}(F) = \inf_{(\lambda_{p})_{p\geq 1}} \sum_{p=1}^\infty g(\diam{\lambda_{p}}),
\]
where the infimum is taken over all sequences $(\lambda_{p})_{p\geq 1}$ in $\Lambda_{g}\cup\{\emptyset\}$ such that $F\subseteq\bigcup_{p}\lambda_{p}$. In addition, for $\overline{g},g\in\jauge_{1}$, let us write $\overline{g}\prec g$ if $\overline{g}/g$ monotonically tends to infinity at zero. We can now give the definition of the classes $\grint^g(\R)$. Recall that a $G_{\delta}$-set is one that may be expressed as a countable intersection of open sets.

\begin{df}[sets with large intersection in $\R$]
For any gauge function $g\in\jauge_{1}$, the class $\grint^g(\R)$ of sets with large intersection in $\R$ with respect to $g$ is the collection of all $G_{\delta}$-subsets $F$ of $\R$ such that $\netm^{\overline{g}}_{\infty}(F\cap U)=\netm^{\overline{g}}_{\infty}(U)$ for every $\overline{g}\in\jauge_{1}$ enjoying $\overline{g}\prec g$ and every open set $U$.
\end{df}

The classes $\grint^g(\T)$ are then defined in the natural way using the classes $\grint^g(\R)$ and the canonical surjection $\phi$ from $\R$ to $\T$.

\begin{df}[sets with large intersection in $\T$]
For any gauge function $g\in\jauge_{1}$, the class $\grint^g(\T)$ of sets with large intersection in $\T$ with respect to $g$ is the collection of all subsets $F$ of $\T$ such that $\phi^{-1}(F)\in\grint^g(\R)$.
\end{df}

The results of \cite{Durand:2007uq} show that the classes $\grint^g(\T)$ of sets with large intersection in the torus enjoy the following remarkable properties.

\begin{thm}\label{grintstable}
For any gauge function $g\in\jauge_{1}$,\begin{enumerate}
\item the class $\grint^g(\T)$ is closed under countable intersections;
\item every set $F\in\grint^g(\T)$ enjoys $\hau^{\overline{g}}(F\cap V)=\infty$ for every $\overline{g}\in\jauge_{1}$ with $\overline{g}\prec g$ and every nonempty open set $V$, and in particular
\[
\dim F\geq s_{g}=\sup\{ s\in (0,1) \:|\: \Id^s\prec g \};
\]
\item every $G_{\delta}$-subset of $\T$ with full Lebesgue measure belongs to $\grint^g(\T)$.
\end{enumerate}
\end{thm}

As previously announced, the sets $\widetilde E_{h}$ defined by (\ref{deftildeEhhmm}) belong to certain classes $\grint^g(\T)$ of sets with large intersection. More precisely, Proposition \ref{mindimEhhmm1} in Section \ref{hmmmultifrac} yields the following result.

\begin{prp}
Let us assume that $\widetilde h$ is finite. Then, with probability one, for all $h\in [\underline{h},\min(\widetilde h,\overline{h}))$, the set $\widetilde E_{h}$ belongs to the class $\grint^{\Id^{h/\widetilde h}}(\T)$.
\end{prp}

Together with Theorem \ref{grintstable}, this result implies that with probability one, for every $h\in [\underline{h},\min(\widetilde h,\overline{h}))$, the set $\widetilde E_{h}$ has infinite Hausdorff measure for every gauge function $g\in\jauge_{1}$ with $g\prec\Id^{h/\widetilde h}$. This property comes into play in the proof of Theorem \ref{multifrac}, because it enables to obtain a sharp lower bound on the Hausdorff dimension of the corresponding iso-H\"older set $E_{h}$, see Section \ref{hmmmultifrac}.

\section{Preparatory lemmas}\label{hmmprel}

In this section, we establish several results that are called upon at various points of the proofs of Theorems \ref{bifrac} and \ref{multifrac}.

The H\"older exponent of the random wavelet series $R$ at a given point $x$ of the torus depends on the way large wavelet coefficients are located around $x$. To be specific, let
\[
S=\{ u\in\etiq \:|\: X_{u}=1 \}=\bigcup_{j=0}^\infty S_{j},
\]
where the sets $S_{j}$ are defined by (\ref{defSj}). As shown by (\ref{defClambda}), the vertices in $S$ correspond to the dyadic intervals indexing the large coefficients of the wavelet series $R$. In addition, for every $u\in\etiq$, let
\begin{equation}\label{defxudotxu}
x_{u}=\phi(\dot x_{u}) \qquad\text{with}\qquad \dot x_{u}=\sum_{j=1}^{\gene{u}} u_{j} 2^{-j}
\end{equation}
and, for every real $\alpha>\underline{h}$, let
\begin{equation}\label{defLalpha}
L_{\alpha}=\{ x\in\T \:|\: \dist(x,x_{u})<2^{-\underline{h}\gene{u}/\alpha}\text{ for infinitely many }u\in S \},
\end{equation}
where $\dist$ is the quotient distance on the torus. It is straightforward to check that $\alpha\mapsto L_{\alpha}$ is nondecreasing. The following lemma establishes a connection between the sets $L_{\alpha}$ and the sets $E_{h}$ and $\widetilde E_{h}$ defined by (\ref{defEhhmm}) and (\ref{deftildeEhhmm}) respectively.

\begin{lem}\label{locholdhmm}
\begin{enumerate}
\item\label{locholdhmm1} For every $h\in [0,\underline{h})\cup(\overline{h},\infty]$, the set $E_{h}$ is empty.
\item\label{locholdhmm2} For every $h\in [\underline{h},\overline{h}]$,
\[
\widetilde E_{h}=\bigcap_{h<\alpha\leq\overline{h}} L_{\alpha} \qquad\text{and}\qquad E_{h}=\widetilde E_{h}\setminus\bigcup_{\underline{h}<\alpha<h} L_{\alpha}.
\]
\end{enumerate}
\end{lem}

\begin{proof}
Assertion (\ref{locholdhmm1}) is due to the fact that the H\"older exponent of the process $R$ is everywhere between $\underline{h}$ and $\overline{h}$, as shown by (\ref{underhleqhRxleqoverh}). Assertion (\ref{locholdhmm2}) follows from the observation that for any $\alpha\in (\underline{h},\overline{h}]$ and any $x\in\T$,
\[
\left\{\begin{array}{rcl}
x\in L_{\alpha} & \Longrightarrow & h_{R}(x)\leq\alpha \\[2mm]
x\not\in L_{\alpha} & \Longrightarrow & h_{R}(x)\geq\alpha.
\end{array}\right.
\]
Indeed, let us assume that $x\in L_{\alpha}$. Then, there are infinitely many dyadic intervals $\lambda=\phi(2^{-j}(k+[0,1)))\in\Lambda$, with $j\in\N_{0}$ and $k\in\intn{0}{2^j-1}$, such that
\[
C_{\lambda}=2^{-\underline{h}\gene{\lambda}} \qquad\text{and}\qquad \dist(x,x_{\lambda})<2^{-\underline{h}\gene{\lambda}/\alpha},
\]
where $x_{\lambda}=\phi(k 2^{-j})$. Proposition \ref{critponclem} ensures that $h_{R}(x)\leq\alpha$. Conversely, let us assume that $x\not\in L_{\alpha}$. Then, for every dyadic interval $\lambda\in\Lambda$ such that $\gene{\lambda}$ is large enough, if $C_{\lambda}=2^{-\underline{h}\gene{\lambda}}$, then $\dist(x,x_{\lambda})\geq 2^{-\underline{h}\gene{\lambda}/\alpha}$. Thus, whether $C_{\lambda}=2^{-\underline{h}\gene{\lambda}}$ or $C_{\lambda}=2^{-\overline{h}\gene{\lambda}}$, we have
\[
|C_{\lambda}|\leq(2^{-\gene{\lambda}}+\dist(x,x_{\lambda}))^\alpha,
\]
so that $h_{R}(x)\geq\alpha$ thanks to Proposition \ref{critponclem}.
\end{proof}

The proof of the following lemma is modeled on that of Proposition 1 in \cite[Chapter 11]{Kahane:1985gc}.

\begin{lem}\label{locholdhmmbis}
With probability one, for every $\alpha\in (\widetilde h,\infty)$, we have $L_{\alpha}=\T$.
\end{lem}

\begin{proof}
Let us consider a real number $\alpha>\widetilde h$. Moreover, for every $u\in\etiq$, let $B_{u}$ be the open ball of $\T$ with center $x_{u}$ and radius $2^{-\underline{h}\gene{u}/\alpha}$. It is straightforward to establish the following inclusion of events:
\begin{equation}\label{incleven1}
\{ \T \neq L_{\alpha} \} \subseteq \bigcup_{j_{0}=0}^\infty \bigcap_{j=j_{0}}^\infty \left\{ \T \neq \bigcup_{u\in S_{j-1}\cup S_{j}} B_{u} \right\}.
\end{equation}

Then, let $j_{0}$ and $j$ be two integers such that
\begin{equation}\label{jj0large}
j-1\geq j_{0}\geq \frac{1-\log_{2}(2^{\underline{h}/\alpha}-1)}{1-\underline{h}/\alpha}
\end{equation}
and let us assume that $\T$ cannot be written as the union over $u\in S_{j-1}\cup S_{j}$ of the balls $B_{u}$. So, there is an integer $k\in\intn{0}{2^j-1}$ such that the closed ball with center $\phi(k 2^{-j})$ and radius $2^{-j-1}$ is not included in this last union of balls. As a result, the point $\phi(k 2^{-j})$ cannot belong to the union over $u\in S_{j-1}\cup S_{j}$ of the open balls with center $x_{u}$ and radius $2^{-\underline{h}\gene{u}/\alpha}-2^{-j-1}$. Therefore,
\begin{equation}\label{incleven2}
\left\{ \T \neq \bigcup_{u\in S_{j-1}\cup S_{j}} B_{u} \right\} \subseteq \bigcup_{k=0}^{2^{j}-1} \left(\mathcal{A}^{j,k}_{j-1} \cap \mathcal{A}^{j,k}_{j}\right),
\end{equation}
where $\mathcal{A}^{j,k}_{j'}$ denotes, for each $j'\in\{j-1,j\}$, the event corresponding to the fact that the Markov chain $X$ maps to $0$ all the vertices of the set
\[
A^{j,k}_{j'}=\{ u\in\{0,1\}^{j'} \:|\: \dist(x_{u},\phi(k 2^{-j}))\leq 2^{-1-\underline{h}j'/\alpha} \}.
\]

Owing to (\ref{jj0large}), there is a set $A'\subseteq A^{j,k}_{j-1}$ such that $\#A'=\lfloor (2^{(1-\underline{h}/\alpha)j}-3)/2\rfloor$ and such that the two sons of each vertex in $A'$ belong to $A^{j,k}_{j}$. The Markov condition (\ref{condMarkov}) then yields
\[
\prob(\mathcal{A}^{j,k}_{j-1}\cap\mathcal{A}^{j,k}_{j})\leq \prob(\forall u\in A' \quad X_{u}=X_{u0}=X_{u1}=0)\leq \nu_{0,j-1}(\{(0,0)\})^{\# A'}.
\]
It follows that the left-hand side of (\ref{incleven2}) is included in an event of probability at most $v_{j}=2^j\exp(-\lfloor (2^{(1-\underline{h}/\alpha)j}-3)/2\rfloor\eta_{j-1})$. Thus, the event
\[
\bigcap_{j=j_{0}}^\infty \left\{ \T \neq \bigcup_{u\in S_{j-1}\cup S_{j}} B_{u} \right\}
\]
has probability at most $\liminf_{j} v_{j}$, which vanishes because there are infinitely many integers $j\geq j_{0}+1$ such that $\eta_{j-1}\geq 2^{(-1+2\underline{h}/(\alpha+\widetilde h))j}$, owing to the fact that $\alpha>\widetilde h$. Then, (\ref{incleven1}) implies that $L_{\alpha}$ is almost surely equal to the whole torus and the result follows from the fact that $\alpha\mapsto L_{\alpha}$ is nondecreasing.
\end{proof}

Let $\alpha\in (\underline{h},\infty)$. Lemma \ref{lemdecompLphhmm} below gives a useful decomposition of the set $L_{\alpha}$ defined by (\ref{defLalpha}). The first set which comes into play in this decomposition is
\begin{equation}\label{defLalphatilde}
\widetilde L_{\alpha}=\{ x\in\T \:|\: \dist(x,x_{u})<2^{-\underline{h}\gene{u}/\alpha}\text{ for infinitely many }u\in\widetilde S \},
\end{equation}
where
\begin{equation}\label{deftildeS}
\widetilde S=\{u\in\etiq^* \:|\: X_{u}=1\text{ and }X_{\pi(u)}=0\}.
\end{equation}
The set $\widetilde S$ is clearly included in $S$, so that $\widetilde L_{\alpha}$ is included in $L_{\alpha}$. Essentially, a point of $L_{\alpha}$ also belongs to $\widetilde L_{\alpha}$ if infinitely many of the large wavelet coefficients $C_{\lambda}$ of the process $R$ that are close to it have appeared at scale $\gene{\lambda}$, i.e.~are such that $C_{\pi(\lambda)}$ is a small wavelet coefficient, where $\pi(\lambda)$ is the smallest dyadic interval of the torus that strictly contains $\lambda$.

The second set coming into play in the decomposition of $L_{\alpha}$ is a set denoted by $\Theta$ and defined as follows. For every vertex $u$ of the binary tree $\etiq$, let $u\etiq$ denote the set of all words of the form $uw$ with $w\in\etiq$ and let
\[
\tau_{u}=\{ v\in u\etiq \:|\: \forall j\in\intn{\gene{u}}{\gene{v}} \quad X_{v_{1}\ldots v_{j}}=1 \}.
\]
The set $\tau_{u}$ is empty if $X_{u}=0$. Otherwise, $\tau_{u}$ is the largest subtree of $\etiq$ rooted at $u$ and formed of vertices which are mapped to $1$ by the Markov chain $X$. The boundary of $\tau_{u}$ is the set
\[
\partial\tau_{u}=\{ \zeta=(\zeta_{j})_{j\geq 1}\in\{0,1\}^{\N} \:|\: \forall j\geq \gene{u} \quad \zeta_{1}\ldots\zeta_{j}\in\tau_{u} \}.
\]
For every sequence $\zeta=(\zeta_{j})_{j\geq 1}$ in $\{0,1\}$, let
\[
\dot x_{\zeta}=\sum_{j=1}^\infty \zeta_{j} 2^{-j}.
\]
Then, let
\begin{equation}\label{defThetahmm}
\dot\Theta=\bigcup_{u\in\etiq}\bigcup_{\zeta\in\partial\tau_{u}} \{\dot x_{\zeta}\} \qquad\text{and}\qquad \Theta=\phi(\dot\Theta),
\end{equation}
where $\phi$ is the canonical surjection from $\R$ to $\T$. Essentially, a point of the torus belongs to $\Theta$ if it can be obtained as the intersection of a sequence of nested dyadic intervals indexing large wavelet coefficients of the process $R$ that propagate across scales. The following lemma provides the law of the Hausdorff dimension of $\Theta$. We refer to Section \ref{results} for the definitions of the parameters appearing in its statement.

\begin{lem}\label{prprecaptheta}
If $\theta<0$, then $\Theta$ is empty with probability one. If not, then with probability one, $\Theta$ is empty or has Hausdorff dimension $\theta$ and, in addition,
\begin{itemize}
\item if there is a $j_{*}\geq 0$ such that $\nu_{1,j}(\{(0,0)\})=0$ for any $j\geq j_{*}$, then
\[
\prob(\Theta=\emptyset)=\Phi_{j_{*}}(0)\cdot\prod_{j=j_{*}}^\infty (1-\eta_{j})^{2^j};
\]
\item if $\nu_{1,j}(\{(0,0)\})>0$ for infinitely many integers $j\geq 0$ and if $\sum_{j}2^j\eta_{j}<\infty$, then $\prob(\Theta=\emptyset)$ is positive and it is equal to one if and only if
\[
\varsigma_{\underline{j}}=\infty
\qquad\text{or}\qquad
\liminf_{j\to\infty}\prod_{\ell=\underline{j}}^j \gamma_{\ell}=0
\qquad\text{or}\qquad
\left\{\begin{array}{l}
\Phi_{\underline{j}}(0)=1 \\[2mm]
\forall j\geq\underline{j} \quad \eta_{j}=0;
\end{array}\right.
\]
\item if $\nu_{1,j}(\{(0,0)\})>0$ for infinitely many integers $j\geq 0$, if $\sum_{j}2^j\eta_{j}=\infty$ and if $\theta>0$, then $\prob(\Theta=\emptyset)$ is less than one and it is equal to zero if and only if $\sum_{j} 2^j\eta_{j}/\varsigma_{j+1}=\infty$.
\end{itemize}
\end{lem}

\begin{proof}
The lemma is a straightforward consequence of Proposition 4 in \cite{Durand:2007lr}, which provides the law of the Hausdorff dimension of the set $\dot\Theta$, and the observation that the sets $\dot\Theta$ and $\Theta$ have the same Hausdorff dimension.
\end{proof}

The following lemma supplies a precise statement of the aforementioned decomposition of the set $L_{\alpha}$ in terms of the sets $\widetilde L_{\alpha}$ and $\Theta$.

\begin{lem}\label{lemdecompLphhmm}
For every $\alpha\in (\underline{h},\infty)$, we have $L_{\alpha}=\widetilde L_{\alpha}\cup\Theta$.
\end{lem}

\begin{proof}
To begin with, it is easy to check that $\widetilde L_{\alpha}\subseteq L_{\alpha}$. Next, let us consider a point $x$ in $\Theta$. Then, there are a vertex $u\in\etiq$ and a sequence $\zeta=(\zeta_{j})_{j\geq 1}\in\partial\tau_{u}$ such that $x=\phi(\dot x_{\zeta})$. For every integer $j\geq\gene{u}$, we have $\zeta_{1}\ldots\zeta_{j}\in S$ and
\[
\dist(x,x_{\zeta_{1}\ldots\zeta_{j}})\leq 2^{-j}<2^{-\underline{h}j/\alpha}.
\]
The point $x$ thus belongs to $L_{\alpha}$. It follows that $\Theta$ is included in $L_{\alpha}$. Hence, $L_{\alpha}$ contains both $\widetilde L_{\alpha}$ and $\Theta$.

Conversely, let us consider a point $x$ in $L_{\alpha}$ which does not belong to $\widetilde L_{\alpha}$. Then, there is an integer $j_{0}\geq 0$ such that $\dist(x,x_{u})\geq 2^{-\underline{h}\gene{u}/\alpha}$ for every vertex $u\in\widetilde S$ with generation at least $j_{0}$. We may assume that $j_{0}\geq (\log_{2}(2^{\underline{h}/\alpha}-1))/(\underline{h}/\alpha-1)$. Let
\[
S'=\{ u\in S \:|\: \dist(x,x_{u})<2^{-\underline{h}\gene{u}/\alpha} \}
\]
and observe that the set $\widetilde S$ cannot contain any vertex of $S'$ with generation at least $j_{0}$. Since $x\in L_{\alpha}$, there exists a sequence $(v^n)_{n\geq 1}$ in $S'$ such that $\gene{v^n}$ is increasing. A standard diagonal argument leads to a sequence $\zeta=(\zeta_{j})_{j\geq 1}$ in $\{0,1\}$ such that for every $j\geq 1$, there are infinitely many integers $n\geq 1$ enjoying $\zeta_{1}\ldots\zeta_{j}=v^n_{1}\ldots v^n_{j}$. Let $u=\zeta_{1}\ldots\zeta_{j_{0}}$ and let us consider two integers $j\geq j_{0}$ and $n$ satisfying $\gene{v^n}>j$ and $\zeta_{1}\ldots\zeta_{j}=v^n_{1}\ldots v^n_{j}$. The vertex $v^n$ belongs to $S'$ and its generation is at least $j_{0}$, so that $v^n\in S\setminus\widetilde S$. Hence, $\pi(v^n)$ belongs to $S$. Moreover,
\[
\dist(x,x_{\pi(v^n)})\leq \dist(x,x_{v^n})+\dist(x_{v^n},x_{\pi(v^n)})<2^{-\underline{h}\gene{v^n}/\alpha}+2^{-\gene{v^n}}\leq 2^{-\underline{h}\gene{\pi(v^n)}/\alpha},
\]
which ensures that $\pi(v^n)\in S'$. By repeating this procedure $\gene{v^n}-j$ times, one can prove that $\zeta_{1}\ldots\zeta_{j}=v^n_{1}\ldots v^n_{j}\in S'$. In particular, $X_{\zeta_{1}\ldots\zeta_{j}}=1$ for every integer $j\geq j_{0}$, so that $\zeta\in\partial\tau_{u}$. Furthermore, for any $j\geq j_{0}$,
\[
\dist(x,\phi(\dot x_{\zeta}))\leq\dist(x,x_{\zeta_{1}\ldots\zeta_{j}})+\dist(x_{\zeta_{1}\ldots\zeta_{j}},\phi(\dot x_{\zeta}))\leq 2^{-\underline{h}j/\alpha}+\sum_{j'=j+1}^\infty \zeta_{j'} 2^{-j'}.
\]
Letting $j\to\infty$, we obtain $x=x_{\zeta}$. The point $x$ thus belongs to $\Theta$.
\end{proof}

\section{Proof of Theorem \ref{bifrac}}

In order to establish Theorem \ref{bifrac}, let us assume that $\sum_{j} 2^j \eta_{j}<\infty$ and $\theta<1$. For any integer $j\geq 1$, let
\[
\widetilde S_{j}=\widetilde S\cap\{0,1\}^j,
\]
where $\widetilde S$ is the set defined by (\ref{deftildeS}). The Markov condition (\ref{condMarkov}) implies that
\begin{equation}\label{majcardtildeSj}\begin{split}
\esp\bigl[\#\widetilde S_{j} \:\bigl|\: \trib_{j-1}\bigr] &=\sum_{u\in\{0,1\}^{j-1}\atop X_{u}=0}\esp\left[\esp[X_{u0}+X_{u1}\:|\:\tribbis_{u}]\:\bigl|\: \trib_{j-1}\right] \\
&\leq 2^{j-1} \int_{\{0,1\}^2} (x_{0}+x_{1})\,\nu_{0,j-1}(\dd x) \leq 2^j\eta_{j-1},
\end{split}\end{equation}
where $\trib_{j-1}$ is the $\sigma$-field generated by the variables $X_{u}$, for $u\in\etiq$ such that $\gene{u}\leq j-1$. It follows that the set $\widetilde S_{j}$ is nonempty with probability at most $2^j \eta_{j-1}$. As the sum $\sum_{j} 2^j \eta_{j}$ converges, the Borel-Cantelli lemma ensures that, with probability one, there are at most finitely many integers $j\geq 1$ such that $\widetilde S_{j}\neq\emptyset$. Consequently, the set $\widetilde S$ is almost surely finite. So, with probability one, for every real $\alpha>\underline{h}$, the set $\widetilde L_{\alpha}$ given by (\ref{defLalphatilde}) is empty. By virtue of Lemmas \ref{locholdhmm} and \ref{lemdecompLphhmm}, it follows that with probability one, for all $h\in [0,\infty]$,
\[
\begin{cases}
E_{h}=\Theta & \text{if }h=\underline{h} \\[2mm]
E_{h}=\T\setminus\Theta & \text{if }h=\overline{h} \\[2mm]
E_{h}=\emptyset & \text{else}.
\end{cases}
\]
Theorem \ref{bifrac} is then a direct consequence of Lemma \ref{prprecaptheta}.

\section{Proof of Theorem \ref{multifrac}}\label{hmmmultifrac}

In order to prove Theorem \ref{multifrac}, let us assume that $\sum_{j} 2^j \eta_{j}=\infty$ and $\theta<1$. To begin with, observe that
\begin{equation}\label{psEhvidehmm}
\text{a.s.} \quad \forall h\in [0,\underline{h})\cup(\min(\widetilde h,\overline{h}),\infty] \qquad E_{h}=\emptyset,
\end{equation}
owing to Lemmas \ref{locholdhmm} and \ref{locholdhmmbis}. Thus, we may now restrict our attention to the case in which $h$ is between $\underline{h}$ and $\min(\widetilde h,\overline{h})$.

The next result gives an upper bound of the Hausdorff dimension of the iso-H\"older set $E_{h}$ for any $h\in [\underline{h},\min(\widetilde h,\overline{h}))$.

\begin{prp}\label{majdimEhhmm}
With probability one,
\[
\dim E_{\underline{h}}\leq\max(\underline{h}/\widetilde h,\dim\Theta) \qquad\text{and}\qquad \forall h\in (\underline{h},\min(\widetilde h,\overline{h})) \quad \dim E_{h}\leq h/\widetilde h.
\]
\end{prp}

\begin{proof}
Owing to (\ref{majcardtildeSj}), the expectation of $\#\widetilde S_{j}$ is at most $2^j \eta_{j-1}$ for any integer $j\geq 1$. Markov's inequality then implies that $\#\widetilde S_{j}$ is greater than $2^j \eta_{j-1} j^2$ with probability at most $1/j^2$. By virtue of the Borel-Cantelli lemma, it follows that with probability one, $\#\widetilde S_{j}$ is bounded by $2^j \eta_{j-1} j^2$ for all $j$ large enough. Thus,
\begin{equation}\label{majtildecardSjkappa}
\text{a.s.} \quad \exists \widetilde\kappa\geq 1 \quad \forall j\geq 1 \qquad \#\widetilde S_{j}\leq\widetilde\kappa\,2^j \eta_{j-1} j^2.
\end{equation}
Let us assume that the event of probability one on which (\ref{majtildecardSjkappa}) holds occurs and let $h\in [\underline{h},\min(\widetilde h,\overline{h}))$ and $s>h/\widetilde h$. For $\alpha\in (h,s\widetilde h)$ and $\eps>0$, the set $\widetilde L_{\alpha}$ defined by (\ref{defLalphatilde}) is covered by the open balls with center $x_{u}$ and radius $2^{-\underline{h}\gene{u}/\alpha}$, for $u\in\widetilde S$ such that $2^{-\underline{h}\gene{u}/\alpha}\leq \eps/2$. Therefore,
\[
\hau^{\Id^s}_{\eps}(\widetilde L_{\alpha})\leq\sum_{j\in\N \atop 2^{-\underline{h}j/\alpha}\leq \eps/2} \#\widetilde S_{j}\cdot(2^{1-\underline{h}j/\alpha})^s\leq\widetilde\kappa\sum_{j\in\N \atop 2^{-\underline{h}j/\alpha}\leq \eps/2} 2^j \eta_{j-1} j^2(2^{1-\underline{h}j/\alpha})^s.
\]
Since $\alpha/s<\widetilde h$, this last series converges so that the right-hand side tends to zero as $\eps\to 0$. Hence, the Hausdorff $\Id^s$-measure of the set $\widetilde L_{\alpha}$ vanishes. It follows that with probability one, for all $h\in [\underline{h},\min(\widetilde h,\overline{h}))$, the Hausdorff dimension of $\bigcap_{h<\alpha\leq\overline{h}} \widetilde L_{\alpha}$ is at most $h/\widetilde h$. The result then follows from the fact that
\[
E_{\underline{h}}=\Theta\cup\bigcap_{\underline{h}<\alpha\leq\overline{h}} \widetilde L_{\alpha} \qquad\text{and}\qquad \forall h\in (\underline{h},\min(\widetilde h,\overline{h})) \quad E_{h}\subseteq\bigcap_{h<\alpha\leq\overline{h}} \widetilde L_{\alpha},
\]
owing to Lemmas \ref{locholdhmm} and \ref{lemdecompLphhmm}.
\end{proof}

Lemma \ref{prprecaptheta} and Proposition \ref{majdimEhhmm}, together with the assumption that $\theta$ is less than one, imply that with probability one, for any $h\in [\underline{h},\min(\widetilde h,\overline{h}))$, the iso-H\"older set $E_{h}$ has Lebesgue measure zero. Moreover, with probability one, this set is empty for every $h\in [0,\underline{h})\cup(\min(\widetilde h,\overline{h}),\infty]$, owing to (\ref{psEhvidehmm}). As a result, with probability one, the set $E_{\min(\widetilde h,\overline{h})}$ has full Lebesgue measure in the torus. In particular, its Hausdorff dimension is equal to one.

In order to give a lower bound on the Hausdorff dimension of the iso-H\"older set $E_{h}$ for every real $h\in [\underline{h},\min(\widetilde h,\overline{h}))$, we shall treat two cases separately: $\widetilde h<\infty$ and $\widetilde h=\infty$. Let us first consider the case in which $\widetilde h$ is finite. The lower bound then follows from the fact that $\widetilde E_{h}$ is a set with large intersection, as shown by the following result. Recall that the classes $\grint^g(\T)$ of sets with large intersection in the torus are defined in Subsection \ref{largeint}.

\begin{prp}\label{mindimEhhmm1}
Let us assume that $\widetilde h<\infty$. Then, with probability one, for every real $h\in [\underline{h},\min(\widetilde h,\overline{h}))$,
\[
\widetilde E_{h}\in\grint^{\Id^{h/\widetilde h}}(\T) \qquad\text{and}\qquad \dim E_{h}\geq h/\widetilde h.
\]
\end{prp}

\begin{proof}
Lemma \ref{locholdhmmbis} shows that, with probability one, for every real $\alpha>\widetilde h$, the set $L_{\alpha}$ defined by (\ref{defLalpha}) is equal to the whole torus. Let us assume that the corresponding event occurs and let $h\in [\underline{h},\min(\widetilde h,\overline{h}))$. For each $\alpha\in (h,\overline{h}]$, we have $\alpha\widetilde h/h>\widetilde h$, so that $L_{\alpha\widetilde h/h}=\T$. Thus, the set
\[
\phi^{-1}(L_{\alpha\widetilde h/h})=\left\{ x\in\R \:\Biggl|\: \begin{array}{l}
|x-p-\dot x_{u}|<2^{-\underline{h}h\gene{u}/(\alpha\widetilde h)}\\
\text{for infinitely many }(u,p)\in S\times\Z
\end{array} \right\},
\]
where $\dot x_{u}$ is defined by (\ref{defxudotxu}), has full Lebesgue measure in $\R$. Following the terminology of \cite{Durand:2007uq}, the family $(p+\dot x_{u},2^{-\underline{h}h\gene{u}/(\alpha\widetilde h)})_{(u,p)\in S\times\Z}$ is a {\em homogeneous ubiquitous system} in $\R$ and Theorem~2 in \cite{Durand:2007uq} implies that the set
\[
\phi^{-1}(L_{\alpha})=\left\{ x\in\R \:\Biggl|\: \begin{array}{l}
|x-p-\dot x_{u}|<2^{-\underline{h}\gene{u}/\alpha}\\
\text{for infinitely many }(u,p)\in S\times\Z
\end{array} \right\}
\]
belongs to the class $\grint^{\Id^{h/\widetilde h}}(\R)$, so that the set $L_{\alpha}$ belongs to the class $\grint^{\Id^{h/\widetilde h}}(\T)$. Furthermore, owing to Lemma \ref{locholdhmm} and the fact that $\alpha\mapsto L_{\alpha}$ is nondecreasing, the set $\widetilde E_{h}$ is equal to the intersection over the integers $n>1/(\overline{h}-h)$ of the sets $L_{h+1/n}$. Each of these sets belongs to the class $\grint^{\Id^{h/\widetilde h}}(\T)$, which is closed under countable intersections thanks to Theorem \ref{grintstable}. Hence, with probability one,
\begin{equation}\label{mindimEhhmm11}
\forall h\in [\underline{h},\min(\widetilde h,\overline{h})) \qquad \widetilde E_{h}\in\grint^{\Id^{h/\widetilde h}}(\T).
\end{equation}

In order to establish the remainder of the proposition, let us begin by observing that $E_{\underline{h}}=\widetilde E_{\underline{h}}\in\grint^{\Id^{\underline{h}/\widetilde h}}(\T)$ with probability one, by virtue of (\ref{mindimEhhmm11}) and Lemma \ref{locholdhmm}. Theorem \ref{grintstable} then implies that the Hausdorff dimension of $E_{\underline{h}}$ is at least $\underline{h}/\widetilde h$ with probability one. Moreover, for every vertex $u\in\etiq$, using a standard diagonal argument, one easily checks that
\[
\bigcup_{\zeta\in\partial\tau_{u}} \{\dot x_{\zeta}\}=\bigcap_{j=\gene{u}}^\infty\bigcup_{v\in\tau_{u} \atop \gene{v}=j} \left( \dot x_{v}+[0,2^{-j}] \right).
\]
This ensures that the set $\dot\Theta$ defined by (\ref{defThetahmm}) is a $F_{\sigma}$-set, i.e.~a set that may be expressed as a countable union of closed sets. Therefore, the set $\Theta=\phi(\dot\Theta)$ is a $F_{\sigma}$-set as well. In addition, this set has Lebesgue measure zero with probability one, because of Lemma \ref{prprecaptheta} and the assumption that $\theta$ is less than one. So, $\T\setminus\Theta$ is almost surely a $G_{\delta}$-subset of $\T$ with full Lebesgue measure. This property, (\ref{majtildecardSjkappa}) and (\ref{mindimEhhmm11}) thus simultaneously hold with probability one. Let us assume that the corresponding event occurs and let $h\in (\underline{h},\min(\widetilde h,\overline{h}))$. Lemmas \ref{locholdhmm} and \ref{lemdecompLphhmm} imply that
\[
E_{h}=(\widetilde E_{h}\setminus\Theta)\setminus\bigcup_{\underline{h}<\alpha<h} \widetilde L_{\alpha}.
\]
In addition, (\ref{mindimEhhmm11}) and Theorem \ref{grintstable} show that the class $\grint^{\Id^{h/\widetilde h}}(\T)$ contains the sets $\widetilde E_{h}$ and $\T\setminus\Theta$ respectively. Since this class is closed under countable intersections, it contains the set $\widetilde E_{h}\setminus\Theta$. Hence, by Theorem \ref{grintstable} again, $\hau^{-\Id^{h/\widetilde h}\cdot\log}(\widetilde E_{h}\setminus\Theta)=\infty$. Furthermore, for $\alpha\in (\underline{h},h)$ and $\eps>0$ small enough, the set $\widetilde L_{\alpha}$ is covered by the open balls with center $x_{u}$ and radius $2^{-\underline{h}\gene{u}/\alpha}$ for $u\in\widetilde S$ such that $2^{-\underline{h}\gene{u}/\alpha}\leq \eps/2$. Thus, owing to (\ref{majtildecardSjkappa}),
\[
\hau^{-\Id^{h/\widetilde h}\cdot\log}_{\eps}(\widetilde L_{\alpha})\leq\widetilde\kappa\sum_{j\in\N \atop 2^{1-\underline{h}j/\alpha}\leq\eps} 2^j \eta_{j-1} j^2(2^{1-\underline{h}j/\alpha})^{h/\widetilde h}\left(\frac{\underline{h}}{\alpha}j-1\right)\log 2.
\]
Since $\alpha\widetilde h/h<\widetilde h$, this last series converges so that the right-hand side tends to zero as $\eps\to 0$. It follows that $\hau^{-\Id^{h/\widetilde h}\cdot\log}(\widetilde L_{\alpha})=0$. The fact that $\alpha\mapsto\widetilde L_{\alpha}$ is nondecreasing implies that
\[
\hau^{-\Id^{h/\widetilde h}\cdot\log}\left(\bigcup_{\underline{h}<\alpha<h} \widetilde L_{\alpha}\right)=0
\]
because the union can actually be taken on a countable subset of $(\underline{h},h)$. Therefore, $\hau^{-\Id^{h/\widetilde h}\cdot\log}(E_{h})=\infty$, so that the Hausdorff dimension of $E_{h}$ is at least $h/\widetilde h$.
\end{proof}

Let us now consider the case in which $\widetilde h$ is infinite. The following result provides a lower bound on the dimension of the iso-H\"older set $E_{h}$, for every $h\in [\underline{h},\overline{h})$.

\begin{prp}\label{mindimEhhmm2}
Let us assume that $\widetilde h=\infty$. Then, with probability one, for every real $h\in [\underline{h},\overline{h})$,
\[
\dim E_{h}\geq 0.
\]
\end{prp}

To prove Proposition \ref{mindimEhhmm2}, we do not use the theory of sets with large intersection. Instead, we follow the main ideas of the proof of Lemma~9 in \cite{Jaffard:1999fg}. Specifically, for any $h\in [\underline{h},\overline{h})$, we obtain a point $y_{h}$ in the set $E_{h}$ as the intersection of a sequence $(I^h_{n})_{n\geq 1}$ of nested closed sets and we show that, with probability one, the construction of this sequence is possible for all $h\in [\underline{h},\overline{h})$. To this end, let us establish the following preparatory result.

\begin{lem}\label{lemconstrucps}
With probability one:\begin{enumerate}
\item\label{lemconstrucps1} For every $u\in\etiq$, the set $u\etiq^*\cap\widetilde S$ is nonempty.
\item\label{lemconstrucps2} There is real $\widetilde\kappa\geq 1$ such that
\[
\forall j\geq 1 \qquad \#\widetilde S_{j} \leq\widetilde\kappa\, 2^j \eta_{j-1} j^2.
\]
\item\label{lemconstrucps3} If $\theta<0$, then the set $\Theta$ is empty. Conversely, if $\theta\geq 0$, then there is a real $\kappa\geq 1$ such that
\[
\forall j\geq 0 \qquad \# S_{j} \leq \kappa\, 2^{(1-\theta)j/2} \prod_{\ell=\underline{j}}^{j-1} \gamma_{\ell}.
\]
\end{enumerate}
\end{lem}

\begin{proof}
To begin with, observe that (\ref{lemconstrucps2}) directly follows from (\ref{majtildecardSjkappa}). So, we only need to prove (\ref{lemconstrucps1}) and (\ref{lemconstrucps3}). In order to prove (\ref{lemconstrucps1}), let $u\in\etiq$. The Markov condition (\ref{condMarkov}) implies that, for any integer $j\geq\gene{u}$,
\[
\prob(\forall v\in u\etiq \quad \gene{v}\leq j \quad\Longrightarrow\quad X_{v}=1)\leq\prod_{\ell=\gene{u}}^{j-1} \nu_{1,\ell}(\{(1,1)\})\leq 2^{\gene{u}-j}\prod_{\ell=\gene{u}}^{j-1}\gamma_{\ell}.
\]
If $\gene{u}\leq\underline{j}$, the right-hand side vanishes for $j$ large enough. Otherwise, since $\theta$ is assumed to be less than one, its limit inferior as $j\to\infty$ vanishes. Thus,
\begin{equation}\label{psexistuxu0}
\text{a.s.} \quad  \exists u' \in u\etiq \qquad X_{u'}=0.
\end{equation}
Furthermore, for any vertex $u'\in\etiq$ and any integer $j>\gene{u'}$, observe that
\begin{equation}\label{incleven4}
\left\{ \forall u''\in u'\etiq \quad X_{u''}=0 \right\} \subseteq \mathcal{B}^{u'}_{j},
\end{equation}
where $\mathcal{B}^{u'}_{j}$ is the event corresponding to the fact that $X_{u''}=0$ for any vertex $u''$ in $u'\etiq$ with generation at most $j$. Due to the Markov condition (\ref{condMarkov}), the conditional probability of the event $\mathcal{B}^{u'}_{j}$ conditionally on the $\sigma$-field generated by the variables $X_{v}$, for $v\in\etiq$ such that $\gene{v}<j$, is equal to $\ind_{\mathcal{B}^{u'}_{j-1}}(1-\eta_{j-1})^{2^{j-\gene{u'}-1}}$, so that the probability of $\mathcal{B}^{u'}_{j}$ is $(1-\eta_{j-1})^{2^{j-\gene{u'}-1}}$ times that of $\mathcal{B}^{u'}_{j-1}$. Arguing by induction, one then readily verifies that
\[
\prob(\mathcal{B}^{u'}_{j})=\prob(\mathcal{B}^{u'}_{\gene{u'}})\prod_{\ell=\gene{u'}+1}^j \left(1-\eta_{\ell-1}\right)^{2^{\ell-\gene{u'}-1}}.
\]
Since $\sum_{\ell} 2^\ell \eta_{\ell}=\infty$, the preceding product tends to zero as $j\to\infty$. By (\ref{incleven4}) and the fact that $\etiq$ is countable, it follows that with probability one, for all $u'\in\etiq$,
\[
X_{u'}=0 \quad \quad\Longrightarrow\quad \quad \exists u''\in u'\etiq^* \qquad X_{u''}=1.
\]
Thanks to (\ref{psexistuxu0}), the set $u\etiq$ almost surely contains a vertex $u'$ enjoying $X_{u'}=0$ and, because of the last assertion, the set $u'\etiq^*$ almost surely contains a vertex $u''$ with $X_{u''}=1$. We may assume that the generation of the vertex $u''$ is minimal. The set $u\etiq^*\cap\widetilde S$, containing $u''$, is thus nonempty. The fact that $\etiq$ is countable finally leads to (\ref{lemconstrucps1}).

Let us prove (\ref{lemconstrucps3}). If $\theta<0$, then Lemma \ref{prprecaptheta} ensures that the set $\Theta$ is almost surely empty. Conversely, let us assume that $\theta\geq 0$. Owing to the Markov condition (\ref{condMarkov}), for any integer $j\geq 1$, the conditional expectation of $\# S_{j}$ conditionally on the $\sigma$-field generated by the variables $X_{u}$ for $u\in\etiq$ with $\gene{u}<j$ is at most $2^j \eta_{j-1}+\gamma_{j-1}\,\# S_{j-1}$. Arguing by induction on $j$, one can establish that
\[
\forall j\geq 0 \qquad \esp[ \# S_{j} ] \leq \sum_{k=-1}^{j-1} 2^{k+1}\eta_{k}\prod_{\ell=k+1}^{j-1} \gamma_{\ell},
\]
with the convention that $\eta_{-1}=1$. As $\gamma_{\underline{j}-1}$ vanishes, it follows that
\[
\forall j\geq\underline{j} \qquad \esp[ \# S_{j} ] \leq \left( \prod_{\ell=\underline{j}}^{j-1} \gamma_{\ell} \right) \sum_{k=\underline{j}-1}^{j-1} \frac{ 2^{k+1} \eta_{k} }{ \prod\limits_{\ell=\underline{j}}^k \gamma_{\ell} }.
\]
Moreover, the fact that $\widetilde h=\infty$ and $\theta\in [0,1)$ implies that for all $k$ large enough,
\[
\eta_{k} \leq 2^{-(7+\theta)k/8} \qquad\text{and}\qquad \prod_{\ell=\underline{j}}^k \gamma_{\ell} \geq 2^{(-1+\theta)(k+1)/8}.
\]
As a result, for some real $c>0$ and every integer $j\geq\underline{j}$, Markov's inequality yields
\[
\prob\left( \# S_{j} > 2^{(1-\theta)j/2 } \prod_{\ell=\underline{j}}^{j-1} \gamma_{\ell} \right) \leq \frac{ \esp\left[ \# S_{j} \right] }{ 2^{(1-\theta)j/2 } \prod\limits_{\ell=\underline{j}}^{j-1} \gamma_{\ell} } \leq c \, 2^{-(1-\theta)j/4}.
\]
We conclude using the Borel-Cantelli lemma.
\end{proof}

From now on, we assume that the event on which the statement of Lemma \ref{lemconstrucps} holds occurs. For any $h\in [\underline{h},\overline{h})$, let us build recursively a sequence $(I^h_{n})_{n\geq 1}$ of nested closed subsets of the torus which lead to a point of the set $E_{h}$. For this purpose, let
\begin{equation}\label{defetab}
\rho^h_{j} = j 2^{\underline{h} j/h} \sum_{j'=j+1}^\infty 2^{(1-\underline{h}/h) j'} \eta_{j'-1} (j')^2
\end{equation}
for any integer $j\geq 0$. Since $\sum_{j} 2^j \eta_{j}=\infty$ and $\widetilde h=\infty$, it is easy to check that $(\rho^h_{j})_{j\geq 0}$ is a sequence of positive reals which enjoys
\begin{equation}\label{etabjcroitlent}
\forall \eps>0 \qquad \rho^h_{j}={\rm o}(2^{\eps j}) \quad\text{as}\quad j\to\infty.
\end{equation}
Moreover, for any vertex $u\in\etiq$, let $B^h_{u}$ be the open ball with center $x_{u}$ and radius $2^{-\underline{h}\gene{u}/h}$ and let $\frac{3}{2} B^{\underline{h}}_{u}$ be the open ball with center $x_{u}$ and radius $3\cdot 2^{-\gene{u}-1}$.

Together with the sequence $(I^h_{n})_{n\geq 1}$ of nested closed sets, we build a nondecreasing sequence $(j^h_{n})_{n\geq 0}$ of nonnegative integers. The construction of the set $I^h_{1}$ and the integers $j^h_{0}$ and $j^h_{1}$ depends on whether or not $\theta$ is negative.

$\bullet$~{\em Step $1$, if $\theta<0$.} Let us build the set $I^h_{1}$ and the integers $j^h_{0}$ and $j^h_{1}$. As $\widetilde h=\infty$, the series $\sum_{j} 2^{(1-\underline{h}/h) j} \eta_{j-1} j^2$ converges, so there is an integer $j^h_{0}\geq 4\widetilde\kappa$ such that
\begin{equation}\label{majsummajcardfjtilde}
\forall j\geq j^h_{0} \qquad 2\widetilde\kappa \sum_{j'=j^h_{0}}^j 2^{(1-\underline{h}/h) j'} \eta_{j'-1} (j')^2 \leq \frac{1}{4},
\end{equation}
where $\widetilde\kappa$ is given by Lemma \ref{lemconstrucps}(\ref{lemconstrucps2}). Furthermore, Lemma \ref{lemconstrucps}(\ref{lemconstrucps1}) shows that, for $j\geq j^h_{0}$ large enough, the set $\widetilde S_{j^h_{0}}\cup\ldots\cup\widetilde S_{j}$ is nonempty. In addition, there is at least one connected component, denoted by $\mathcal{I}_{j}$, of the complement in the torus of the balls $B^h_{u}$, for $u\in\widetilde S_{j^h_{0}}\cup\ldots\cup\widetilde S_{j}$, which has Lebesgue measure at least
\[
\frac{ 1 - \sum\limits_{j'=j^h_{0}}^j \#\widetilde S_{j'} \, 2^{1-\underline{h} j'/h} }{ \sum\limits_{j'=j^h_{0}}^j \#\widetilde S_{j'} } \geq \frac{ 3 }{ 4\widetilde\kappa \sum\limits_{j'=j^h_{0}}^j 2^{j'} \eta_{j'-1} (j')^2 }.
\]
Note that this inequality follows from Lemma \ref{lemconstrucps}(\ref{lemconstrucps2}) and (\ref{majsummajcardfjtilde}). The component $\mathcal{I}_{j}$ can contain the image under the canonical surjection $\phi$ of some closed subinterval of $\R$ with length $\rho^h_{j} 2^{-\underline{h}j/h}$ if
\[
\rho^h_{j} \leq \frac{ 3\cdot 2^{\underline{h}j/h} }{ 4\widetilde\kappa \sum\limits_{j'=j^h_{0}}^j 2^{j'} \eta_{j'-1} (j')^2 }.
\]
As $\widetilde h=\infty$, the right-hand side tends to infinity exponentially fast as $j\to\infty$. Thus, (\ref{etabjcroitlent}) implies that the preceding inequality holds for $j$ large enough. Let $j^h_{1}$ be the smallest integer such that $\widetilde S_{j^h_{0}}\cup\ldots\cup\widetilde S_{j^h_{1}}$ is nonempty and the inequality holds. Remark that, next to $\mathcal{I}_{j^h_{1}}$, there is a ball $B^h_{u}$ with $u\in\widetilde S_{j^h_{0}}\cup\ldots\cup\widetilde S_{j^h_{1}}$. Hence, $\mathcal{I}_{j^h_{1}}$ contains a set, denoted by $I^h_{1}$, of the form
\[
\phi( x_{u} + 2^{-\underline{h}\gene{u}/h} + [ 0 , \rho^h_{j^h_{1}} 2^{-\underline{h} j^h_{1}/h} ]) \quad\text{or}\quad \phi( x_{u} - 2^{-\underline{h}\gene{u}/h} + [ -\rho^h_{j^h_{1}} 2^{-\underline{h} j^h_{1}/h} , 0 ]),
\]
with $u\in\widetilde S_{j^h_{0}}\cup\ldots\cup\widetilde S_{j^h_{1}}$. Clearly, the intersection of the sets $B^h_{v}$ and $I^h_{1}$ is empty for every vertex $v\in\widetilde S_{j^h_{0}}\cup\ldots\cup\widetilde S_{j^h_{1}}$.

$\bullet$~{\em Step~$1$, if $\theta\geq 0$.} Since $\theta\in [0,1)$, there is an infinite subset $\mathcal{J}$ of $\N$ such that $\gamma_{\underline{j}}\ldots \gamma_{j-2} \leq 2^{(5\theta+1)(j-1)/6}$ for any integer $j\in\mathcal{J}$. Together with Lemma \ref{lemconstrucps}(\ref{lemconstrucps3}), this ensures that $\# S_{j-1}$ is bounded by $\kappa \, 2^{(2+\theta)(j-1)/3}$ for any $j\in\mathcal{J}$. Let $j^h_{0}$ denote an integer in $\mathcal{J}$ which is greater than $4\widetilde\kappa$ and is large enough to ensure both (\ref{majsummajcardfjtilde}) and
\begin{equation}\label{majcardfjdelta}
3\kappa \, 2^{-(1-\theta)( j^h_{0} - 1 )/3 } \leq \frac{1}{4}.
\end{equation}
By virtue of Lemma \ref{lemconstrucps}(\ref{lemconstrucps1}), for any $j\geq j^h_{0}$ large enough, the set $\widetilde S_{j^h_{0}}\cup\ldots\cup\widetilde S_{j}$ is nonempty. In addition, there is at least one connected component, denoted by $\mathcal{I}_{j}$, of the complement in $\T$ of the balls $\frac{3}{2} B^{\underline{h}}_{u}$, for $u\in S_{j^h_{0}-1}$, and the balls $B^h_{u}$, for $u\in\widetilde S_{j^h_{0}}\cup\ldots\cup\widetilde S_{j}$, which has Lebesgue measure at least
\[
\frac{ 1 - 3\#S_{j^h_{0}-1}\, 2^{-( j^h_{0}-1 )} - \sum\limits_{j'=j^h_{0}}^j \#\widetilde S_{j'} \, 2^{1-\underline{h} j'/h} }{ \#S_{j^h_{0}-1} + \sum\limits_{j'=j^h_{0}}^j \#\widetilde S_{j'} } \geq \frac{ 1/2 }{\kappa \, 2^{\frac{2+\theta}{3}(j^h_{0}-1)} + \widetilde\kappa \sum\limits_{j'=j^h_{0}}^j 2^{j'} \eta_{j'-1} (j')^2}.
\]
This last inequality follows from Lemma \ref{lemconstrucps}(\ref{lemconstrucps2}-\ref{lemconstrucps3}), along with (\ref{majsummajcardfjtilde}) and (\ref{majcardfjdelta}). The component $\mathcal{I}_{j}$ can contain the image under $\phi$ of some closed subinterval of $\R$ with length $\rho^h_{j} 2^{-\underline{h}j/h}$ if
\[
\rho^h_{j} \leq \frac{ 2^{-1+\underline{h}j/h} }{\kappa \, 2^{(2+\theta)(j^h_{0}-1)/3} + \widetilde\kappa \sum\limits_{j'=j^h_{0}}^j 2^{j'} \eta_{j'-1} (j')^2 }.
\]
This inequality holds for $j$ large enough, because of (\ref{etabjcroitlent}) and the fact that $\widetilde h$ is infinite. Let $j^h_{1}$ denote the smallest integer for which $\widetilde S_{j^h_{0}}\cup\ldots\cup\widetilde S_{j^h_{1}}$ is nonempty and the inequality holds. Next to $\mathcal{I}_{j^h_{1}}$, there is a ball $\frac{3}{2} B^{\underline{h}}_{u}$ with $u\in S_{j^h_{0}-1}$ or a ball $B^h_{u}$ with $u\in\widetilde S_{j^h_{0}}\cup\ldots\cup\widetilde S_{j^h_{1}}$. Thus, $\mathcal{I}_{j^h_{1}}$ contains a set, denoted by $I^h_{1}$, of the form
\[
\phi(x_{u} + 3\cdot 2^{-j^h_{0}} + [0,\rho^h_{j^h_{1}} 2^{-\underline{h}j^h_{1}/h}]) \quad\text{or}\quad \phi(x_{u} - 3\cdot 2^{-j^h_{0}} +[-\rho^h_{j^h_{1}} 2^{-\underline{h} j^h_{1}/h},0]),
\]
with $u\in S_{j^h_{0}-1}$, or of the form
\[
\phi(x_{u} + 2^{-\underline{h}\gene{u}/h} + [0,\rho^h_{j^h_{1}} 2^{-\underline{h} j^h_{1}/h}]) \quad\text{or}\quad \phi(x_{u} - 2^{-\underline{h}\gene{u}/h} + [-\rho^h_{j^h_{1}} 2^{-\underline{h} j^h_{1}/h},0]),
\]
with $u\in\widetilde S_{j^h_{0}}\cup\ldots\cup\widetilde S_{j^h_{1}}$. Note that $B^h_{v} \cap I^h_{1}=\emptyset$ for every vertex $v\in\widetilde S_{j^h_{0}}\cup\ldots\cup\widetilde S_{j^h_{1}}$ and that $\frac{3}{2} B^{\underline{h}}_{v} \cap I^h_{1}=\emptyset$ for every vertex $v\in S_{j^h_{0}-1}$.

$\bullet$~{\em Step $n+1$ for $n\geq 1$.} Steps $1$ to $n$ have supplied the sets $I^h_{1}\supseteq\ldots\supseteq I^h_{n}$ and the integers $j^h_{0}\leq\ldots\leq j^h_{n}$. Let us build the set $I^h_{n+1}$ and the integer $j^h_{n+1}$. Because of (\ref{defetab}) and the fact that $j^h_{n}\geq j^h_{0}\geq 4\widetilde\kappa$, we have
\begin{equation}\label{majrestsummajcardfjtilde}
\forall j\geq j^h_{n}+1 \qquad 2\widetilde\kappa\sum_{j'=j^h_{n}+1}^j 2^{(1-\underline{h}/h) j'} \eta_{j'-1} (j')^2 \leq \frac{\rho^h_{j^h_{n}} 2^{-\underline{h} j^h_{n}/h}}{2}.
\end{equation}
Let us consider a vertex $v\in\etiq$ enjoying $\gene{v}\geq j^h_{n}+1$ and $\overline{\lambda_{v}}\subseteq I^h_{n}$. Lemma \ref{lemconstrucps}(\ref{lemconstrucps1}) ensures that there is a vertex $v'\in v\etiq^*\cap\widetilde S$. As a result, the set $B^h_{v'}\cap I^h_{n}$, containing the point $x_{v'}$, is nonempty and the set $\widetilde S_{j^h_{n}+1}\cup\ldots\cup\widetilde S_{j}$ is nonempty for $j$ large enough. In addition, there is at least one connected component, denoted by $\mathcal{I}_{j}$, of the complement in $I^h_{n}$ of the balls $B^h_{u}$, for $u\in\widetilde S_{j^h_{n}+1}\cup\ldots\cup\widetilde S_{j}$, which has Lebesgue measure at least
\[
\frac{ \rho^h_{j^h_{n}} 2^{-\underline{h} j^h_{n}/h} - \sum\limits_{j'=j^h_{n}+1}^j \#\widetilde S_{j'} \, 2^{1-\underline{h} j'/h} }{ 1+\sum\limits_{j'=j^h_{n}+1}^j \#\widetilde S_{j'} }\geq \frac{ \rho^h_{j^h_{n}} 2^{-\underline{h} j^h_{n}/h} }{ 2\left( 1+\widetilde\kappa \sum\limits_{j'=j^h_{n}+1}^j 2^{j'} \eta_{j'-1} (j')^2 \right) }.
\]
The inequality follows from Lemma \ref{lemconstrucps}(\ref{lemconstrucps2}) and (\ref{majrestsummajcardfjtilde}). The component $\mathcal{I}_{j}$ can contain the image under $\phi$ of some closed subinterval of $\R$ with length $\rho^h_{j} 2^{-\underline{h} j/h}$ if
\[
\rho^h_{j} \leq \frac{ \rho^h_{j^h_{n}} 2^{\underline{h} (j - j^h_{n})/h} }{ 2\left( 1+\widetilde\kappa \sum\limits_{j'=j^h_{n}+1}^j 2^{j'} \eta_{j'-1} (j')^2 \right) },
\]
which holds for $j$ large enough because of (\ref{etabjcroitlent}) and the fact that $\widetilde h$ is infinite. Let $j^h_{n+1}$ be the smallest integer such that this inequality holds and such that there exists a vertex $v'\in\widetilde S_{j^h_{n}+1}\cup\ldots\cup\widetilde S_{j^h_{n+1}}$ with $x_{v'}\in I^h_{n}$. Observe that, next to $\mathcal{I}_{j^h_{n+1}}$, there is a ball $B^h_{u}$ with $u\in\widetilde S_{j^h_{n}+1}\cup\ldots\cup\widetilde S_{j^h_{n+1}}$. As a consequence, $\mathcal{I}_{j^h_{n+1}}$ contains a set, denoted by $I^h_{n+1}$, of the form
\[
\phi(x_{u} + 2^{-\underline{h}\gene{u}/h} + [0,\rho^h_{j^h_{n+1}} 2^{-\underline{h} j^h_{n+1}/h}]) \quad\text{or}\quad \phi(x_{u} - 2^{-\underline{h}\gene{u}/h} +[-\rho^h_{j^h_{n+1}} 2^{-\underline{h} j^h_{n+1}/h},0])
\]
with $u\in\widetilde S_{j^h_{n}+1}\cup\ldots\cup\widetilde S_{j^h_{n+1}}$. Observe that the intersection of the sets $B^h_{v}$ and $I^h_{n+1}$ is empty for every vertex $v\in\widetilde S_{j^h_{n}+1}\cup\ldots\cup\widetilde S_{j^h_{n+1}}$.

The sets $I^h_{0},I^h_{1},\ldots$ given by the preceding procedure form a decreasing sequence of closed subsets of the torus. Moreover, the diameter of each set $I^h_{n}$ is at most $\rho^h_{j^h_{n}} 2^{-\underline{h} j^h_{n}/h}$, which tends to zero as $n\to\infty$ by virtue of (\ref{etabjcroitlent}). Thus, the intersection over $n\geq 1$ of the sets $I^h_{n}$ is a singleton $\{y_{h}\}$.

\begin{lem}
The point $y_{h}$ belongs to the iso-H\"older set $E_{h}$.
\end{lem}

\begin{proof}
Let $\alpha\in (h,\overline{h}]$. Owing to (\ref{etabjcroitlent}), there exists an integer $n_{0}\geq 2$ such that $j^h_{n_{0}-1}\geq(\log_{2} 3)/(\underline{h}/h-\underline{h}/\alpha)$ and $\rho^h_{j^h_{n}} 2^{-\underline{h} j^h_{n}/h}\leq 2^{-\underline{h} j^h_{n}/\alpha}/3$ for every integer $n\geq n_{0}$. Let us consider an integer $n\geq n_{0}$. The point $y_{h}$ belongs to $I^h_{n}$, so there exists a vertex $u^n\in\widetilde S_{j^h_{n-1}+1}\cup\ldots\cup\widetilde S_{j^h_{n}}$ for which
\[
\dist(y_{h},x_{u^n}) \leq 2^{-\underline{h}\gene{u^n}/h}+\rho^h_{j^h_{n}} 2^{-\underline{h} j^h_{n}/h} < 2^{-\underline{h}\gene{u^n}/\alpha}
\]
The last inequality follows from the fact that $n\geq n_{0}$ and $\gene{u^n}\leq j^h_{n}$. As a result, the point $y_{h}$ belongs to the set $L_{\alpha}$ defined by (\ref{defLalpha}). By Lemma \ref{locholdhmm}, this point thus belongs to the set $\widetilde E_{h}$.

This proves the lemma in the case where $h=\underline{h}$. We may therefore assume that $h>\underline{h}$. Lemma \ref{locholdhmm}, together with the fact that $\alpha\mapsto L_{\alpha}$ is nondecreasing, shows that it suffices to establish that $y_{h}\not\in L_{h}$.

Let us assume that $\theta<0$. The point $y_{h}$ belongs to $I^h_{1}$, so it cannot belong to any ball $B^h_{u}$ for $u\in\widetilde S_{j^h_{0}}\cup\ldots\cup\widetilde S_{j^h_{1}}$. Moreover, for any integer $n\geq 1$, the point $y_{h}$ belongs to $I^h_{n+1}$, so it cannot belong to any ball $B^h_{u}$ for $u\in\widetilde S_{j^h_{n}+1}\cup\ldots\cup\widetilde S_{j^h_{n+1}}$. It follows that $y_{h}$ does not belong to any ball $B^h_{u}$ with $u\in\widetilde S$ and $\gene{u}\geq j^h_{0}$. Hence, $y_{h}$ does not belong to the set $\widetilde L_{h}$ defined by (\ref{defLalphatilde}). Furthermore, Lemma \ref{lemdecompLphhmm} shows that $L_{h}=\widetilde L_{h}\cup\Theta$ and Lemma \ref{lemconstrucps}(\ref{lemconstrucps3}) implies that $\Theta$ is empty. It follows that $y_{h}\not\in L_{h}$.

Let us assume that $\theta\geq 0$. In this case, the point $y_{h}$ does not belong to any ball $B^h_{u}$ with $u\in\widetilde S$ and $\gene{u}\geq j^h_{0}$ and does not belong to any closed dyadic interval $\overline{\lambda_{u}}$ with $u\in S_{j^h_{0}-1}$. Therefore, the point $y_{h}$ cannot belong to the set $\widetilde L_{h}$. It cannot belong to the set $\Theta$ either. Otherwise, there would exist a vertex $u\in\widetilde S$ and a sequence $\zeta=(\zeta_{j})_{j\geq 1}\in\partial\tau_{u}$ enjoying $y_{h}=\phi(\dot x_{\zeta})$. If $\gene{u}\leq j^h_{0}-1$, then the vertex $\zeta_{1}\ldots\zeta_{j^h_{0}-1}$ would belong to $S_{j^h_{0}-1}$ and index a closed dyadic interval containing $y_{h}$. If $\gene{u}\geq j^h_{0}$, then the point $y_{h}$ would belong to the set $\overline{\lambda_{u}}$ and thus to the ball $B^h_{u}$, since $h>\underline{h}$. In both cases, we would end up with a contradiction. Hence, $y_{h}$ does not belong to $\Theta$. Lemma \ref{lemdecompLphhmm} finally ensures that $y_{h}\not\in L_{h}$.
\end{proof}

We have established that, with probability one, for any $h\in [\underline{h},\overline{h})$, it is possible to build a point $y_{h}$ in the set $E_{h}$. Proposition \ref{mindimEhhmm2} is thus proven.

Using Propositions \ref{majdimEhhmm}, \ref{mindimEhhmm1} and \ref{mindimEhhmm2} together with the fact that $E_{\underline{h}}\supseteq\Theta$ by Lemmas \ref{locholdhmm} and \ref{lemdecompLphhmm}, we finally obtain the following result.

\begin{prp}
With probability one,
\[
\dim E_{\underline{h}}=\max(\underline{h}/\widetilde h,\dim\Theta) \qquad\text{and}\qquad \forall h\in (\underline{h},\min(\widetilde h,\overline{h})) \quad \dim E_{h}=h/\widetilde h.
\]
\end{prp}

Theorem \ref{multifrac} is then an immediate consequence of this result and Lemma \ref{prprecaptheta}.

\section{Proof of Proposition \ref{prpexposchmm}}\label{oscsing}

In order to prove Proposition \ref{prpexposchmm}, let $h\in [\underline{h},\overline{h}]$ and let $x\in E_{h}$. Let us first assume that $h<\overline{h}$ and let us consider a real number $\beta>h/\underline{h}-1$. Owing to Lemma \ref{locholdhmm}, for any integer $n\geq 1$ such that $h+1/n<(\beta+1)\underline{h}$, there exists a dyadic interval $\lambda_{n}\in\Lambda$ enjoying $\gene{\lambda_{n}}\geq n$, $C_{\lambda_{n}}=2^{-\underline{h}\gene{\lambda_{n}}}$ and $\dist(x,x_{\lambda_{n}})<2^{-\underline{h}\gene{\lambda_{n}}/(h+1/n)}$. Note that these intervals $\lambda_{n}$ are such that $d(x,x_{\lambda_{n}})^{1+\beta}\leq 2^{-\gene{\lambda_{n}}}$. Proposition~3 in \cite{Arneodo:1998ai} then shows that $\beta_{R}(x)\leq\beta$. This inequality holds for any $\beta>h/\underline{h}-1$. Thus, $\beta_{R}(x)\leq h/\underline{h}-1$.

Conversely, let us consider a real number $\beta>\beta_{R}(x)$. Owing to Proposition~3 in \cite{Arneodo:1998ai}, there exists a sequence $(\lambda_{n})_{n\geq 1}$ of dyadic intervals of the torus such that $d(x,x_{\lambda_{n}})^{1+\beta}\leq 2^{-\gene{\lambda_{n}}}$ for all $n\geq 1$,
\[
2^{-\gene{\lambda_{n}}}+\dist(x,x_{\lambda_{n}})\xrightarrow[n\to\infty]{} 0 \qquad\text{and}\qquad \frac{\log|C_{\lambda_{n}}|}{\log(2^{-\gene{\lambda_{n}}}+\dist(x,x_{\lambda_{n}}))}\xrightarrow[n\to\infty]{} h.
\]
As $h<\overline{h}$, for infinitely many integers $n\geq 1$, we have
\[
\frac{\log_{2}|C_{\lambda_{n}}|}{-\gene{\lambda_{n}}}\leq\frac{\log|C_{\lambda_{n}}|}{\log(2^{-\gene{\lambda_{n}}}+d(x,x_{\lambda_{n}}))}<\overline{h},
\]
so that $X_{u_{\lambda_{n}}}=1$. Thanks to Lemma \ref{locholdhmm}, it follows that $h\leq (1+\beta)\underline{h}$. Letting $\beta\to\beta_{R}(x)$, we obtain $\beta_{R}(x)\geq h/\underline{h}-1$.

Let us now suppose that $h=\overline{h}<\infty$ and consider a real $\beta>0$. As $h>\underline{h}$, the point $x$ does not belong to $L_{\underline{h}}$ by virtue of Lemma \ref{locholdhmm}. Hence, there exists an integer $j_{0}\geq 0$ such that $X_{u}=0$ for any vertex $u\in\etiq$ enjoying $\gene{u}\geq j_{0}$ and $d(x,x_{u})<2^{-\gene{u}}$. For any integer $j\geq j_{0}$, there is a vertex $u^j\in\etiq$ satisfying $\gene{u^j}=j$ and $d(x,x_{u^j})<2^{-j}$. Observe that
\[
2^{-\gene{\lambda_{u^j}}}+\dist(x,x_{\lambda_{u^j}})\xrightarrow[j\to\infty]{} 0 \qquad\text{and}\qquad \frac{\log|C_{\lambda_{u^j}}|}{\log(2^{-\gene{\lambda_{u^j}}}+\dist(x,x_{\lambda_{u^j}}))}\xrightarrow[j\to\infty]{} h.
\]
Proposition~3 in \cite{Arneodo:1998ai} then ensures that $\beta\geq\beta_{R}(x)$. We conclude by letting $\beta\to 0$.

\end{document}